\documentclass{amsart}
\usepackage{graphicx,amssymb}
\usepackage{hyperref}

\newcommand{\mapdown}[1]%
{\Big\downarrow\rlap{$\vcenter{\hbox{$\scriptstyle#1$}}$}}

\newcommand{\dom}{\mathop{\mathrm{Dom}}}
\newcommand{\img}[1]{\mathop{\mathrm{IMG}}\left(#1\right)}
\newcommand{\arr}{\longrightarrow}
\newcommand{\xs}{X^*}
\newcommand{\symm}[1][X]{\mathfrak{S}\left(#1\right)}
\newcommand{\autxs}{\mathop{\mathrm{Aut}}(\xs)}
\newcommand{\M}{\mathcal{M}}
\newcommand{\Ms}{\mathcal{C}}
\newcommand{\mc}{\mathcal{G_\Ms}}
\newcommand{\F}{\mathfrak{F}}
\newcommand\Gx{{\widetilde G}}
\newcommand{\C}{\mathbb{C}}
\newcommand{\N}{\mathbb{N}}
\newcommand{\nuke}{\mathcal N}
\newcommand{\CS}{\widehat{\C}}
\newcommand{\Sph}{\mathbb S^2}
\newcommand{\Z}{\mathbb{Z}}
\newcommand{\Einfty}{{\mathcal E_\infty}}
\newcommand\opsi{{\overline\psi}}
\newcommand{\ophi}{{\overline\phi}}
\newbox\pairbox
\def\pair<#1,#2>{{\mathsurround=0pt
    \setbox\pairbox\hbox{$\left\langle#1,\;#2\right\rangle$}
    \left\langle\kern-0.35\ht\pairbox
      \copy\pairbox\kern-0.35\ht\pairbox\right\rangle}}

\newtheorem{theorem}{Theorem}[section]
\newtheorem{proposition}[theorem]{Proposition}
\newtheorem{corollary}[theorem]{Corollary}
\newtheorem{lemma}[theorem]{Lemma}

\theoremstyle{definition}
\newtheorem{defi}{Definition}[section]

\begin{document}
\title{Thurston equivalence of topological polynomials}
\author{Laurent Bartholdi}
\author{Volodymyr Nekrashevych}
\date{Preliminary version, last modified by LB, 20060522; compiled \today}
\thanks{L.B.\ and V.N.\ gratefully acknowledge support from Beijing
  University and EPFL, respectively}
\begin{abstract}
  We answer Hubbard's question on determining the Thurston equivalence
  class of ``twisted rabbits'', i.e.\ composita of the ``rabbit''
  polynomial with $n$th powers of the Dehn twists about its ears. The
  answer is expressed in terms of the 4-adic expansion of $n$.

  We also answer the equivalent question for the other two families of
  degree-$2$ topological polynomials with three post-critical points.

  In the process, we rephrase the questions in group-theoretical
  language, in terms of wreath recursions.
\end{abstract}
\maketitle
\tableofcontents

\section{Introduction}
Consider the ``rabbit'' polynomial $f_R(z)\approx z^2+(-0.1226 +
0.7449i)$, whose critical point $0$ is on a periodic orbit of
length $3$. Up to affine transformations, there are exactly two
other polynomials with same action on post-critical points (with
the same \emph{ramification graph}), called the ``corabbit''
$f_C\approx z^2+(-0.1226 - 0.7449i)$ and the ``airplane''
$f_R\approx z^2-1.7549$. Furthermore, by a result of W.\ Thurston (see
below), every branched covering with same ramification graph is
equivalent to precisely one of $f_R,f_C,f_A$.

Consider now a Dehn twist $T$ of $\C$ around the two non-critical
values of the $f_R$-orbit of $0$. The map $T^mf_R$ is again a branched
covering, and it has the same ramification graph as $f_R$; therefore
it is equivalent (i.e., conjugate up to homotopies) to one of
$f_R,f_C,f_A$.  Which one?

This question was asked by J.\ Hubbard, see~\cite{pilgrim:endom}.
The answer, as we shall show (Theorem~\ref{th:4adic}), is the
following:
\begin{quote}
  Write $m$ in base $4$, as $m=\sum_{i=0}^\infty m_i4^i$ with
  $m_i\in\{0,1,2,3\}$ and almost all $m_i=0$ if $m$ is non-negative,
  and almost all $m_i=3$ if $m$ is negative. If one of the $m_i$
  is $1$ or $2$, then $T^mf_R$ is equivalent to $f_A$. Otherwise,
  it is equivalent to $f_R$ if $m$ is non-negative, and to $f_C$ if $m$
  is negative.
\end{quote}
For example, $T^{-1}f_R$ and $T^{-4}f_R$ are equivalent to $f_C$,
since $-1=\ldots 333$ and $-4=\ldots 3330$.

Consider now the polynomial $f_i=z^2+i$, whose finite critical point
has orbit $0\mapsto i\mapsto i-1\mapsto -i\mapsto i-1$.  A branched
covering with that ramification graph is either equivalent to $f_i$,
or to $f_{-i}$, or is not equivalent to any rational map (it is
\emph{obstructed}). A.\ Douady and J.\ Hubbard ask
(\cite[page~293]{DH:Thurston}) to determine, as a function of a Dehn
twist $D$, when $f_i\cdot D$ is obstructed, and if not, whether it is
equivalent to $f_i$ or to $f_{-i}$. The answer (see
Theorem~\ref{th:mod5}) depends on the image of $D$ in a finite group
of order $100$.

Thurston's theorem does not tell us when two obstructed maps
are equivalent. We may however also answer that question:
Corollary~\ref{cor:mod5} shows that there are infinitely many
inequivalent obstructed maps with the same ramification graph as
$f_i$, and gives an algorithm to determine equivalence among
obstructed maps.

The first construction of infinitely many non-equivalent Thurston maps
with the same ramification graph was presented
in~\cite[Proposition~2.12]{census}.

\subsection{Thurston's theorem}
Consider more generally a branched covering $f$ of $\Sph$, with set of
critical points $C_f$, and let $P_f$ be its \emph{post-critical set}:
$P_f=\bigcup_{n\ge1}f^{\circ n}(C_f)$.  Let us suppose that $P_f$ is
finite, in which case $f$ is called \emph{post-critically finite}, or
a \emph{Thurston map}. Two Thurston maps $f,g$ are \emph{equivalent}
if there exist orientation-preserving homeomorphisms
$\phi_0,\phi_1:(\Sph,P_f)\arr(\Sph,P_g)$ such that $\phi_0$ and
$\phi_1$ are isotopic relative to $P_f$, and $\phi_0f=g\phi_1$.
Recall also that a branched covering is a \emph{topological
  polynomial} if $f^{-1}(\infty)=\{\infty\}$.

A \emph{multicurve} is a system
$\Gamma=\{\gamma_1,\dots,\gamma_n\}$ of simple, closed, disjoint,
non-homotopic, non-peripheral curves on $\Sph\setminus P_f$. Here
a curve is called \emph{peripheral} if one of the two parts into
which it divides the sphere contains less than two post-critical
points. The multicurve $\Gamma$ is \emph{stable} if for all
$\gamma\in\Gamma$, all non-peripheral elements of $f^{-1}(\gamma)$
are homotopic to elements of $\Gamma$. There is then an induced
map
\[f_\Gamma:\mathbb{R}^\Gamma\arr\mathbb{R}^\Gamma,\quad\gamma_i\mapsto\sum_{\delta\in
  f^{-1}(\gamma_i)}\frac{[\delta]}{\deg f|_\delta},\] where $[\delta]$
is, if it exists, the element of $\Gamma$ homotopic to $\delta$, and
is $0$ otherwise.

\begin{theorem}[Thurston's criterion;
  see~\cite{DH:Thurston}]\label{th:thurston}
  A Thurston map $f$ with hyperbolic orbifold\footnote{For a
    definition of a hyperbolic orbifold see~\cite{DH:Thurston}.}  is
  equivalent to a rational function if and only if the spectral radius
  of $f_\Gamma$ is $<1$ for all stable multicurves $\Gamma$.

 In that case, the rational function equivalent to $f$ is unique
up to conjugation by a M\"obius transformation.
\end{theorem}

A stable multicurve $\Gamma$ is an \emph{obstruction} if the spectral radius
of $f_\Gamma$ is greater than or equal to $1$. Therefore, Thurston's Theorem says that
a Thurston map is equivalent to a rational map if and only if there are no
obstructions.

If the Thurston map $f$ is a topological polynomial,
then the structure of obstructions is better understood.  Namely,
every Thurston obstruction of a topological polynomial contains a
\emph{Levy cycle} (see~\cite{bielhubb:polynom}), i.e., a multicurve
$\{\gamma_0, \gamma_1, \ldots, \gamma_{n-1}\}$ such that the only
non-peripheral component $\widetilde\gamma_{i-1}$ of
$f^{-1}(\gamma_i)$ is homotopic to $\gamma_{i-1}$ and the maps
$f:\widetilde\gamma_{i-1}\arr\gamma_i$ are of degree $1$ for every $i$
(here the indices $i$ are considered modulo $n$).

Thurston's theorem does not, in principle, provide an
algorithmic answer to the question when a Thurston map $f$ is
equivalent to a rational function; nor does it construct the
rational function. Many attempts were made to that end,
notably~\cite{pilgrim:endom} and~\cite{kam:thurston}. The present
paper may be seen as another step in this direction.

\subsection{Sketch of the method}
Given a post-critically finite branched covering $f$, we associate
with it a finitely generated group acting faithfully on a rooted tree.
Its action is given recursively by self-similar tree isometries. This
group, the \emph{iterated monodromy group} of $f$, is an invariant for
the branched covering up to equivalence.\footnote{It is a complete
  invariant for combinatorial equivalence, if equipped with additional
  algebraic data --- see~\cite[Theorem~6.5.2]{nek:book} and
  Proposition~\ref{pr:homotopy}.} In favourable
(``contracting'') cases, the ``nucleus'' of its action is a
finite-state automaton characterising the group. This gives an
effective method to solve J.  Hubbard's question for any given
$m$.

Let $P$ be the post-critical set of $f$. The (pure) mapping class
group $\mc$ of $\Ms=\C\setminus P$ acts on the branched coverings with
post-critical set $P$ by pre- and post-composition.  The action of
$\mc$ can be pushed through the recursion to give a self-map $\opsi$
(almost a homomorphism) of $\mc$ such that $f\cdot g$ and $f\cdot
\opsi(g)$ are combinatorially equivalent.  The map $\opsi$ is
contracting in the case of the ``rabbit'' polynomial: there is a
finite set $L\subset\mc$ such that for any $g\in\mc$ we have
$\opsi^n(g)\in L$ for some $n$.

It remains to compute the iterated monodromy group of $f\cdot g$
for all $g\in L$ to obtain a general answer to Douady's and
Hubbard's questions.

\subsection{Outline of the paper}
We explain in Section~\ref{se:img} the construction of iterated
monodromy groups, the fundamental notion of ``contracting''
actions, and nuclei. We apply that construction to the study of
Thurston equivalence and study the post- and pre-composition
actions of the mapping class group on the covering in Section~\ref{se:pcf}.

We consider in more detail in Section~\ref{se:tr} the ``rabbit''
polynomial, for which we obtain the recursions as described above,
and we specialise in Subsection~\ref{se:m} to Dehn twists around
the rabbit's ears, obtaining an answer to Hubbard's question.

We re-express the solution in Section~\ref{se:moduli} in more
classical terms of iterations of the pull-back map on Teichm\"uller
and moduli spaces. This gives another approach to the problem, which
is however less algorithmical.

We finally consider in Sections~\ref{se:3} and~\ref{s:quater} the
other two classes of degree-two topological polynomials with three
finite post-critical points.

Section~\ref{se:3} deals with the polynomials whose ramification graph
is the same as that of $z^2+i$ (period $2$, preperiod $1$). In particular,
we classify the obstructed examples up to combinatorial equivalence.
The last section considers the case of period $1$ and preperiod $2$.


\subsection{A remark on notation} We compose in most cases
transformations as if they acted on the right: in a product
$f\cdot g$ the transformation $f$ acts before $g$. In particular,
if $g_1$ and $g_2$ are elements of a group, then
$g_1^{g_2}=g_2^{-1}g_1g_2$. However, we sometimes have to also
consider left actions. Therefore, if we write $f\circ g$, then we
mean left composition, in which $g$ acts before $f$.

\section{Iterated monodromy groups}\label{se:img}
We give here an overview of techniques and results in self-similar
actions and iterated monodromy groups. More details and proofs can
be found in~\cite{nek:book}.

\subsection{Partial self-coverings and monodromy action}
A covering is a continuous surjective map $f: \Ms_1\arr \Ms$ of
topological spaces, with the property that for every $z\in \Ms$ there
exists an open neighborhood $\mathcal U$ of $z$ such that
$f^{-1}(\mathcal U)$ is a disjoint union of open sets each of which is
mapped homeomorphically onto $\mathcal U$ by $f$. We say that the
covering $f$ is \emph{$d$-fold} if there exists $d\in\N$ such
that $|f^{-1}(z)|=d$ for all $z\in \Ms$.

A \emph{partial self-covering} is a covering map $f:\Ms_1\arr\Ms$
such that $\Ms_1$ is an open subspace of a path connected and
locally path connected space $\Ms$.

Let $f:\Ms_1\arr\Ms$ be a $d$-fold partial self-covering. Then the
$n$th iteration $f^n:\Ms_n\arr\Ms$ of $f$ is a $d^n$-fold partial
self-covering with domain $\Ms_n$, usually smaller than $\Ms_1$.

Choose a basepoint $t\in\Ms$. Then the fundamental group $\pi_1(\Ms,
t)$ acts naturally by monodromy on each of the sets $f^{-n}(t)$. The
action it induces on the disjoint union $\bigsqcup_{n\ge 0}f^{-n}(t)$
is called the \emph{iterated monodromy action}.

The set $\bigsqcup_{n\ge 0}f^{-n}(t)$ has a natural structure of a
rooted $d$-regular tree. The root of this tree is the point $t\in
f^{-0}(t)=\{t\}$ and each vertex $z\in f^{-n}(t)$, for $n\ge 1$,
is connected to the vertex $f(z)\in f^{-(n-1)}(t)$.

The \emph{iterated monodromy group} of the partial self-covering
$f$, denoted $\img{f}$, is the quotient of the fundamental group
$\pi_1(\Ms, t)$ by the kernel of the iterated monodromy action.
This kernel, by definition, consists of the loops
$\gamma\in\pi_1(\Ms, t)$ such that for every $n$ all
$f^n$-preimages of $\gamma$ are also loops.

In other words, the iterated monodromy group of $f$ is the group
of all automorphisms of the tree $\bigsqcup_{n\ge 0} f^{-n}(t)$
which are induced by elements of $\pi_1(\Ms, t)$.

The iterated monodromy group can be effectively computed using the
following recursive formula. Choose an alphabet $X$ of $d$
letters. Then every rooted $d$-regular tree is isomorphic to the
tree $\xs$ of finite words over the alphabet $X$. In this tree the
empty word is the root and every word $v\in\xs$ is connected by
edges to all the words of the form $vx$, where $x\in X$ is a
letter.

We will use throughout the paper the following notation. If $\gamma$
is a path in $\Ms$ and $z$ is an $f$-preimage of the startpoint of
$\gamma$, then
\begin{equation}\label{eq:notation}f^{-1}(\gamma)[z]\end{equation}
denotes the unique $f$-preimage of $\gamma$ starting at $z$.

Choose a bijection $\Lambda:X\arr f^{-1}(t)$ and a path $\ell_x$
in $\Ms$ starting at $t$ and ending in $\Lambda(x)$ for every
$x\in X$. Extend the bijection $\Lambda:X\arr f^{-1}(t)$ to
an isomorphism of the rooted trees $\Lambda:\xs\arr\bigsqcup_{n\ge
0}f^{-n}(t)$ inductively by the condition that
\[\Lambda(xv)\text{ be the end of the path }f^{-|v|}(\ell_x)[\Lambda(v)].\]

It is not hard to prove that $\Lambda$ is a well defined
isomorphism of the rooted trees
(see~\cite[Proposition~5.2.1]{nek:book}). We then have the
following recursive description of the iterated monodromy action
(see~\cite[Proposition~5.2.2]{nek:book}).

\begin{proposition}
  \label{pr:imgrecursion} Define an action of $\pi_1(\Ms, t)$ on $\xs$
  by conjugating the iterated monodromy action by $\Lambda$.  Then for
  all $\gamma\in\pi_1(\Ms, t)$, $v\in\xs$ and $x\in X$ we have
  \[(xv)^\gamma=y\left(v^{(\ell_x\gamma_x\ell_y^{-1})}\right),\] where
  $\Lambda(y)$ is the end of the path
  $\gamma_x=f^{-1}(\gamma)[\Lambda(x)]$.
\end{proposition}

We multiply the paths in the natural order: in a product
$\gamma_1\gamma_2$ the path $\gamma_1$ is followed before the path
$\gamma_2$. Therefore $\pi_1(\Ms, t)$ and $\img{f}$ act from the
right on the tree of preimages and on $\xs$. (Note that this is
different from the convention in~\cite{nek:book}.)

\subsection{Wreath recursions} The recursive formula in
Proposition~\ref{pr:imgrecursion} can be interpreted either as a
\emph{wreath recursion}, or as the description of an automaton.

A wreath recursion is a homomorphism $\Phi$ from a group $G$ to
the wreath product $G\wr\symm$, where $\symm$ is the symmetric
group on $X$. Recall that the wreath product $G\wr\symm$ is, by
definition, the semidirect product $G^X\rtimes\symm$, where
$\symm$ acts on the direct power $G^X$ by permutation of the
coordinates. If $X=\{1,\dots,d\}$, we write the elements of the
wreath product in the form $\pair<g_1, g_2, \ldots, g_d>\pi$,
where $\pair<g_1, g_2, \ldots, g_d>$ is an element of the direct
power $G^X$ and $\pi$ is an element of the symmetric group
$\symm$. The elements of the wreath product are multiplied
according to the rule
\[\pair<g_1, g_2, \ldots g_d>\rho\pair<h_1, h_2, \ldots, h_d>\tau=
\pair<g_1h_{1^{\rho}}, g_2h_{2^{\rho}}, \ldots,
g_dh_{d^{\rho}}>\rho\tau.\]

We denote by $g|_x$ the $x$th coordinate of $\Phi(g)$, for $g\in
G$ and $x\in X$. Inductively, we put
\[g|_{vx}=\left(g|_v\right)|_x\]
for all $v\in\xs$ and $x\in X$.

We let $G$ act on $X$ by postcomposing $\Phi$ with the natural
homomorphism $G\wr\symm\arr\symm$. We extend this action to the
\emph{associated action of $G$ on $\xs$}, defined by the recursion
\[(xv)^g=x^g\left(v^{g|_x}\right).\]

\noindent The associated action and restrictions satisfy the following
relations
\begin{equation}\label{eq:restrictions}(gh)|_v=g|_vh|_{v^g},\quad
g|_{vw}=g|_v|_w\end{equation} for all $g, h\in G$ and $v,
w\in\xs$.

If $X=\{1, 2, \ldots, d\}$, then Proposition~\ref{pr:imgrecursion}
can be expressed in terms of wreath recursions in the following
way:
\begin{proposition}
\label{pr:imgwr} The action of $\pi_1(\Ms, t)$ on $\xs$ is the
action associated with the wreath recursion $\Phi:\pi_1(\Ms,
t)\arr\pi_1(\Ms, t)\wr\symm$ given by
\[\Phi(\gamma)=\pair<\ell_1\gamma_1\ell_{k_1}^{-1}, \ell_2\gamma_2\ell_{k_2}^{-1},
\ldots, \ell_d\gamma_d\ell_{k_d}^{-1}>\rho,\] where
$\gamma_i=f^{-1}(\gamma)[\Lambda(i)]$, and $\Lambda(k_i)$ is the
endpoint of $\gamma_i$, and $\rho$ is the permutation $i\mapsto
k_i$.
\end{proposition}

\subsection{Virtual endomorphisms}\label{ss:virend}
A wreath recursion $\Phi:G\arr G\wr\symm$ can be constructed using
the \emph{associated virtual endomorphism}. A virtual endomorphism
$\phi:G\dashrightarrow G$ is a homomorphism $\dom\phi\arr G$ from
a subgroup of finite index into $G$. The subgroup $\dom\phi$ is
called the \emph{domain} of the virtual endomorphism $\phi$.

If $\Phi:G\arr G\wr\symm$ is a wreath recursion and $x_0\in X$ is
a letter, then the domain of the associated virtual endomorphism
$\phi=\phi_{x_0}$ is the stabilizer of the letter $x_0$ (with
respect to the action of $G$ on $X$); and the virtual endomorphism
is defined as the restriction
\[\phi(g)=g|_{x_0}.\]

It follows from~\eqref{eq:restrictions} that $\phi:\dom\phi\arr G$
is a homomorphism.

Suppose that the action of $G$ on $X=\{1, 2, \ldots, d\}$ is
transitive, choose $x_0=1$, and choose some $r_i\in G$, for all $i\in
X$, such that $1^{r_i}=i$. Write $h_i=r_i|_1$. Then the wreath
recursion can be reconstructed by the formula
\[\Phi(g)=\pair<h_1^{-1}\phi(r_1gr_{k_1}^{-1})h_{k_1}, h_2^{-1}\phi(r_2gr_{k_2}^{-1})h_{k_2},
\ldots, h_d^{-1}\phi(r_dgr_{k_d}^{-1})h_{k_d}>\rho,\] where
$h_i=r_i|_1$, and
 $\rho$ is the permutation $i\mapsto k_i$, and the indices $k_i$
are uniquely defined by the condition
$r_igr_{k_i}^{-1}\in\dom\phi$.

If we change the elements $h_i$, or if we change $\{r_i\}$ to
another left coset transversal, then we change $\Phi$ to its
post-composition by an inner automorphism of $G\wr\symm$, and
therefore conjugate the associated actions of $G$ on $X^*$ by an
automorphism of $X^*$. More precisely (see~\cite[Sections~2.3
and~2.5]{nek:book}), there will exist an automorphism $\Delta$ of
the $d$-regular rooted tree, conjugating $G$ into the new action;
and $\Delta$ will satisfy a recursion of the form
\[\Phi(\Delta)=\pair<g_1\Delta,\ldots,g_d\Delta>\pi,\]
where $\pair<g_1, \ldots, g_d>\pi$ is the element defining the
inner automorphism of $G\wr\symm$.

\subsection{Contraction} A wreath recursion $\Phi:G\arr G\wr\symm$ is
\emph{contracting} if there exists a finite set $\nuke\subset G$ such
that for every $g\in G$ there exists $n_0\in\N$ with $g|_v\in\nuke$
for all words $v\in\xs$ of length greater than $n_0$. This property
ensures that many calculations regarding arbitrary elements of $G$ can
be reduced, via the wreath recursion, to considerations on a finite
set. For example, the ``word problem'' (determining if a given product
of $N$ generators is trivial) can be solved in polynomial time in
contracting groups.

If $G$ is generated by a finite symmetric set $S=S^{-1}$, then a
subset $\nuke\subset G$ satisfies the above condition if and only
if $1\in\nuke$ and there exists $n_0\in\N$ such that
\[(gs)|_v\in\nuke\] for all $g\in\nuke, s\in S$ and words $v\in\xs$ of
length greater than $n_0$.

The smallest set $\nuke$ satisfying the condition of these definitions
is called the \emph{nucleus} of the contracting action.

If $\phi:G\dashrightarrow G$ is a virtual endomorphism of a
finitely generated group, then its \emph{spectral radius} is equal
to
\[r(\phi)=\limsup_{n\to\infty} \sqrt[n]{\limsup_{g\in\dom\phi^n,
l(g)\to\infty}\frac{l\left(\phi^n\left(g\right)\right)}{l(g)}},\]
where $l(g)$ denotes the word length of $g$ with respect to some
fixed generating set of $G$.

\def\0{\cite[Proposition~2.11.11]{nek:book}}
\begin{proposition}[\0]
Let  $\Phi:G\arr G\wr\symm$ be a wreath recursion and let $\phi$
be an associated virtual endomorphism.

If $\Phi$ is contracting, then $r(\phi)<1$. If the action of $G$
on $\xs$ is transitive on every level $X^n$ (in particular, if it
is transitive on $X^1$ and $\phi$ is onto) and $r(\phi)<1$, then
the wreath recursion $\Phi$ is contracting.
\end{proposition}

It is proved in~\cite[Theorem~5.5.3]{nek:book} that if a partial
self-covering $f:\Ms_1\arr\Ms$ is expanding, then the associated
wreath recursion $\Phi_f$ on $\pi_1(\Ms)$, as defined in
Proposition~\ref{pr:imgwr}, is contracting.

\subsection{Automata} It is convenient to describe wreath recursions and
nuclei of contracting wreath recursions in terms of
\emph{automata}.

A subset $A\subset G$ is \emph{state-closed} if for every $g\in A$ and
$x\in X$ we have $g|_x\in A$. It is easy to see that the nucleus of a
contracting wreath recursion is a state-closed set.

If $A$ is a state-closed set, then we interpret it as an automaton,
which, when it is in a state $g\in A$ and it reads a letter $x\in X$
on the input tape, prints the letter $x^g$ on the output tape and goes
to the state $g|_x$. Then the automaton $A$ with initial state $g$
transforms any word $v\in\xs$ to the word $v^g$ and thus describes the
associated action of $G$ on $\xs$.

We draw state-closed sets as graphs (\emph{Moore diagrams}) with
vertex set $A$. The vertex $g\in A$ is marked by its image in
$\symm$, i.e., by the permutation it induces on $X$, and for every
$g\in A$ and $x\in X$ we draw an arrow from $g$ to $g|_x$ labeled by
$x$. Then the graph completely describes the restriction of the
wreath recursion to $A$.

We will always have in our paper $X=\{0, 1\}$; then the symmetric
group $\symm$ consists of two elements: $1$ and $\sigma=(0, 1)$.  If
an element of $G^X$ or $\symm$ is trivial, then we usually do not
write it, so that $\pair<g_0, g_1>=\pair<g_0, g_1>1$ and
$\sigma=\pair<1, 1>\sigma$. The elements of $G\wr\symm$ are either
written $g=\pair<g_0, g_1>$ (they are then called \emph{inactive}), or
$g=\pair<g_0, g_1>\sigma$ (and are then called \emph{active}). All
computations in $G\wr\symm$ are based then on two rules:
\[\pair<g_0, g_1>\cdot\pair<h_0, h_1>=\pair<g_0h_0, g_1h_1>\quad\text{and}
\quad\sigma\cdot\pair<g_0, g_1>=\pair<g_1, g_0>\sigma.\] In
drawing Moore diagrams, we indicate active states by a grey dot
and inactive states by a white dot.

\section{Post-critically finite topological polynomials}\label{se:pcf}
\subsection{Homotopy in terms of wreath recursion}
\label{ss:homotopy}

Let $P_0\subset\C$ be a finite set of complex numbers and suppose
that a map $f_P:P_0\arr P_0$ and a point $c_0\in P_0$ are given
such that $P_0=\{f_P^{\circ n}(c_0)\;:\;n\ge 0\}$. Let
$\mathcal{F}$ be the set of all degree-two orientation-preserving
branched coverings $f:\CS\arr\CS$, with critical points $c_0$ and
$\infty$, whose restriction to $P_0$ coincides with $f_P$ and
$f^{-1}(\infty)=\{\infty\}$.

Let us denote by $P$ the set of \emph{post-critical points} of
$f$, i.e., $P=\{f^{\circ n}(c_0)\;:\;n\ge 1\}\cup\{\infty\}$. Note
that $P=P_0\cup\{\infty\}$ if $c_0$ is periodic (i.e., if
$f_P^{\circ n}(c_0)=c_0$ for some positive $n$) and
$P=\left(P_0\setminus\{c_0\}\right)\cup\{\infty\}$ otherwise. The
\emph{ramification graph} of $f$ is the directed graph with vertex
set $P$ and an arrow from $p$ to $f(p)$ for each $p\in P$.

If $f\in\mathcal{F}$, then $\img{f}$ is defined as the iterated
monodromy group of $f:\Ms_1\arr\Ms$, where $\Ms=\CS\setminus P$
and $\Ms_1=f^{-1}\left(\Ms\right)$.

Let us denote by $\F$ the set of homotopy classes (within
$\mathcal F$) of branched coverings $f\in\mathcal{F}$; in other
words, $\F$ is the set of path-connected components of the space
$\mathcal{F}$.

Choose a basepoint $t\in\Ms$, imagined close to infinity. Note
that $\pi_1(\Ms, t)$ is a free group of rank $|P|-1$; we can take
its generators to be loops going around the finite post-critical
points in the positive direction.

Let $a$ be a small simple closed loop in $\CS$ going around
$\infty$ in the negative direction and based at $t$ (on $\C$ it is
a big loop going in the positive direction around all the finite
post-critical points). The loop $a$ divides the sphere $\CS$ into
two parts. One contains $\infty$ and the other contains
$P\setminus\{\infty\}$.

We may assume, after changing the map $f$ to a homotopic one, that
$f(a)=a$. Then $f$ maps $a$ onto itself by a degree-two covering.
We call $a$ the \emph{circle at infinity}. We may also assume that
$f(t)=t$ (again after changing $f$ to a homotopic map). Let us
denote $t$ by $+\infty$ and the $f$-preimage of $t$ different from
$t$ by $-\infty$. The loop $a$ has two $f$-preimages. One is the
subcurve of $a$ starting at $+\infty$ and ending at $-\infty$. The
other is the curve starting at $-\infty$ and ending at $+\infty$.

Let the connecting path $\ell_0$ be the trivial path starting and
ending in the basepoint $+\infty$. Let $\ell_1$ be the $f$-preimage of
$a$ starting at $+\infty$ and ending in $-\infty$ (the \emph{upper
  semicircle at infinity}). Let $\Phi_f:\pi_1(\Ms,
t)\arr\symm\wr\pi_1(\Ms, t)$ be the wreath recursion defined by $f$
and the given choice of connecting paths and let $\Lambda_f:\pi_1(\Ms,
t)\arr\autxs$ be the associated iterated monodromy action on the tree
$\xs$ (see Proposition~\ref{pr:imgwr} and the definition of the
associated action before it).

Note that $\Lambda_f\left(a\right)$ is always equal to the
standard adding machine $\tau=\pair<1,\tau>\sigma$.

\begin{proposition}
\label{pr:homotopy} Two branched coverings $f_0,
f_1\in\mathcal{F}$ are homotopic if and only if there exists
$n\in\Z$ such that for every $\gamma\in\pi_1(\Ms, +\infty)$ we
have
\[\Phi_{f_0}\left(\gamma\right)=\Phi_{f_1}(a^{n}\cdot\gamma\cdot
a^{-n}).\]
\end{proposition}

\begin{proof}
Suppose that the branched coverings $f_0$ and $f_1$ are homotopic.
Let $f_x$ be a homotopy between them, where $x$ varies from $0$ to
$1$. The basepoint $t=+\infty$ has two preimages under $f_x$. The
first preimage $t_0(x)$ draws, as $x$ ranges over $[0,1]$, a path
$g_0$ starting at $+\infty$. The other preimage $t_1(x)$ draws a
path $g_1$ starting at $-\infty$.

Either $t_0(1)=+\infty$ and $t_1(1)=-\infty$, and then $g_0, g_1$
are loops, or $t_0(1)=-\infty$ and $t_1(1)=+\infty$, and then
$g_0$ goes from $+\infty$ to $-\infty$ and $g_1$ goes from
$-\infty$ to $+\infty$. We have chosen $t=+\infty$ to be close to
the point $\infty\in P$ and the homotopies must fix the point
$\infty$, so $t_0(x)$ and $t_1(x)$ remain close to $\infty$.
Consequently, if $g_0, g_1$ are loops, then $g_0=a^n$ and
$g_1=\ell_1a^m\ell_1^{-1}$ for some $n, m\in\Z$, and if $g_0, g_1$
are not loops, then they are of the form $a^n\ell_1$ and
$\overline{\ell_1}a^m$, where $\overline{\ell_1}$ is the lower
semi-circle at infinity (we have $a=\ell_1\overline{\ell_1}$).
Given paths $b,c$, we write $b=c$ here and below to indicate that they are
homotopic.

One of the $f_x$-preimages of the loop $a$ is a path $\ell(x,
-):[0, 1]\arr\Ms$ starting in $t_0(x)$ and ending in $t_1(x)$. We
have $\ell(0, -)=\ell_1$. The path $\ell(1, -)$ is equal to
$\ell_1$ if $g_0$ and $g_1$ are loops, and to the lower semicircle
$\overline{\ell_1}$ otherwise. Note also $\ell(-, 0)=g_0$ and
$\ell(-, 1)=g_1$.

The path $\ell(x, y)$ depends continuously on $x$, and therefore
defines a continuous map $\ell(x, y):[0, 1]\times [0, 1]\arr\Ms$.
It follows that if $g_0$ and $g_1$ are loops, then
$g_0=\ell_1\cdot g_1\cdot\ell_1^{-1}$, and if not, then
$g_0=\ell_1\cdot g_1\cdot\left(\overline{\ell_1}\right)^{-1}$.

Consequently, if $g_0=a^n$ then $g_1=\ell_1^{-1}a^n\ell_1$, and if
$g_0=a^n\ell_1$ then $g_1=\ell_1^{-1}a^n\ell_1
\overline{\ell_1}=\ell_1^{-1}a^{n+1}$.

Take an arbitrary loop $\gamma\in\pi_1(\Ms, t)$, and let
$\gamma_0=f_0^{-1}(\gamma)[+\infty]$ and
$\gamma_1=f_0^{-1}\left(\gamma\right)[-\infty]$ denote their
preimages starting at $+\infty$ and $-\infty$ respectively. The
homotopy from $f_0$ to $f_1$ deforms the paths $\gamma_i$
continuously, giving paths $\gamma_{0,
x}=f_x^{-1}\left(\gamma\right)\left[t_0(x)\right]$ and $\gamma_{1,
x}=f_x^{-1}\left(\gamma\right)\left[t_1(x)\right]$.

Suppose first that $\gamma_0$, $\gamma_1$ and $g_0=a^n,
g_1=\ell_1^{-1}a^n\ell_1$ are loops. Then we get a homotopy
transforming $\gamma_0$ via $t_0\left([0, x]\right)\gamma_{0,
x}t_0\left([0, x]\right)^{-1}$ to the loop
\[g_0\gamma_{0}'g_0^{-1}=a^n\gamma_0'a^{-n},\]
where $\gamma_0'=f_1^{-1}\left(\gamma\right)[+\infty]$.

Similarly, the loop $\ell_1\gamma_1\ell_1^{-1}$ is homotopic to
\[\ell_1g_1\gamma_1'g_1^{-1}\ell_1^{-1}
=a^n\ell_1\gamma_1'\ell_1^{-1}a^{-n},\] where
$\gamma_1'=f_1^{-1}\left(\gamma\right)[-\infty]$.

Since we assumed $\gamma_i$ are loops,
$\Phi_{f_1}(\gamma)=\pair<\gamma_0',\ell_1\gamma_1'\ell_1^{-1}>$.
We then compute
\begin{multline*}\Phi_{f_0}\left(\gamma\right)=
\pair<\gamma_0,\ell_1\gamma_1\ell_1^{-1}>=
\pair<a^{n}\gamma_0'a^{-n},a^n\ell_1\gamma_1'\ell_1^{-1}a^{-n}>\\
=\pair<a^n,a^n>\cdot\Phi_{f_1}(\gamma)\cdot\pair<a^{-n},a^{-n}>=
\Phi_{f_1}\left(a^{2n}\cdot\gamma\cdot
a^{-2n}\right).\end{multline*}

If on the other hand
$\Phi_{f_1}(\gamma)=\pair<\widetilde{\gamma}_0',\widetilde{\gamma}_1'>\sigma$,
then $\Phi_{f_1}(\gamma\cdot a)= \pair<\widetilde{\gamma}_0',
\widetilde{\gamma}_1'>\sigma\pair<1, a>\sigma=
\pair<\widetilde{\gamma}_0'a,\widetilde{\gamma}_1'>$. We then
obtain $\Phi_{f_0}\left(\gamma
a\right)=\Phi_{f_1}\left(a^{2n}\cdot \gamma a\cdot
a^{-2n}\right)$, so again
\[\Phi_{f_0}\left(\gamma\right)=\Phi_{f_1}\left(a^{2n}\cdot\gamma\cdot a^{-2n}\right).\]

Consider now the case $g_0=a^n\ell_1$ and
$g_1=\ell_1^{-1}a^{n+1}$. Take a loop $\gamma\in\pi_1(\Ms,
+\infty)$ such that its $f_1$-preimages are loops $\gamma_0$ and
$\gamma_1$, so that we have
$\Phi_{f_0}\left(\gamma\right)=\pair<\gamma_0,\ell_1\gamma_1\ell_1^{-1}>$.
Let $\gamma_0'$ and $\gamma_1'$ be the $f_1$-preimages of
$\gamma$, so that
$\Phi_{f_1}\left(\gamma\right)=\pair<\gamma_0',\ell_1\gamma_1'\ell_1^{-1}>$.
Then, as before, $\gamma_0$ is homotopic to
$g_0\gamma_1'g_0^{-1}=a^n\ell_1\gamma_1' \ell_1^{-1}a^{-n}$, and
$\ell_1\gamma_1\ell_1^{-1}$ is homotopic to
$\ell_1g_1\gamma_0'g_1^{-1} \ell_1^{-1}
=\ell_1\ell_1^{-1}a^{n+1}\gamma_0'a^{-n-1}\ell_1 \ell_1^{-1}
=a^{n+1}\gamma_0'a^{-n-1}$. We then have
\begin{multline*}\Phi_{f_0}\left(\gamma\right)=\pair<\gamma_0,
\ell_1\gamma_1\ell_1^{-1}>= \pair<a^n\ell_1\gamma_1'
\ell_1^{-1}a^{-n},a^{n+1}\gamma_0'a^{-n-1}>\\
=\pair<a^n,a^{(n+1)}>\sigma\cdot
\pair<\gamma_0',\ell_1\gamma_1'\ell_1^{-1}>\cdot\sigma\pair<a^{-n},a^{-(n+1)}>=
\Phi_{f_1}\left(a^{2n+1}\cdot\gamma\cdot
a^{-2n-1}\right).\end{multline*}

If $\Phi_{f_1}(\gamma)$ is of the form
$\pair<\widetilde{\gamma}_0',\widetilde{\gamma}_1'>\sigma$, we reduce
to the previous case by multiplying by $a$, as before. We also obtain
$\Phi_{f_0}\left(\gamma\right)=\Phi_{f_1}\left(a^{2n+1}\cdot\gamma\cdot
  a^{-2n-1}\right)$.

Suppose now, in order to prove the proposition in the other
direction, that there exists $n\in\Z$ such that
\begin{equation}
\label{eq:phia}
\Phi_{f_0}\left(\gamma\right)=\Phi_{f_1}\left(a^n\cdot\gamma\cdot
a^{-n}\right)\end{equation} holds for all
$\gamma\in\pi_1\left(\Ms, +\infty\right)$.

For $i\in\{0,1\}$, denote by $G_i\le\pi_1(\Ms, +\infty)$
the set of loops whose preimages under $f_i$ are loops. The set
$G_i$ is an index-two subgroup, isomorphic both to
$\pi_1\left(f_i^{-1}\left(\Ms\right), +\infty\right)$ and to
$\pi_1\left(f_i^{-1}\left(\Ms\right), -\infty\right)$, where the
isomorphisms $G_i\arr\pi_1\left(f_i^{-1}\left(\Ms\right),
\pm\infty\right)$ are the maps $\gamma\mapsto
f_i^{-1}(\gamma)[\pm\infty]$. Note that if $\gamma\in G_i$, then
\[\Phi_{f_i}\left(\gamma\right)=\left(f_i^{-1}(\gamma)[+\infty],\quad
\ell_1\cdot
f_{i}\left(\gamma\right)[-\infty]\cdot\ell_1^{-1}\right).\]

Condition~\eqref{eq:phia} implies that $G_0=G_1$ and that for all
$\gamma\in G_0=G_1$:
\begin{eqnarray}\label{eq:vareps}
f_{0}^{-1}\left(\gamma\right)[+\infty] &\text{is
homotopic to}& f_1^{-1}\left(a^{n}\gamma a^{-n}\right)[+\infty]\text{ in }\Ms,\text{ and}\\
\label{eq:vareps2} \ell_1\cdot
f_0^{-1}\left(\gamma\right)[-\infty]\cdot \ell_1^{-1} &\text{is
homotopic to}& \ell_1\cdot f_1^{-1}\left(a^n\gamma
a^{-n}\right)[-\infty]\cdot\ell_1^{-1}.
\end{eqnarray}

The equality $G_0=G_1$ implies that there are homeomorphisms $h_+,
h_-:f_0^{-1}\left(\Ms\right)\arr f_1^{-1}\left(\Ms\right)$ such
that $f_0=f_1\circ h_{\pm}$, where $h_+$ fixes the preimages
$+\infty$, $-\infty$ of the basepoint, and $h_-$ permutes them.

The homeomorphisms $h_\pm$ fix then the points $c_0, c_1, c_2$ and
can be extended in a unique way to homeomorphisms
$\widetilde{h_\pm}:\Ms\arr\Ms$.


Suppose first that $n$ is even. Then
$f_0^{-1}(a^n)[+\infty]=a^{n/2}$ and we have
from~\eqref{eq:vareps}
\[\widetilde{h_+}\left(f_0^{-1}\left(\gamma
\right)[+\infty]\right)=f_1^{-1}\left(\gamma
\right)[+\infty]=f_0^{-1}\left(a^{-n}\gamma
a^n\right)[+\infty]=a^{-n/2}\cdot f_0^{-1}(\gamma)[+\infty]\cdot
a^{n/2},\] where all equalities are homotopies in $\Ms$. We then
have
\[\widetilde{h_+}(\gamma)=a^{-n/2}\gamma a^{n/2}\] for all
$\gamma\in\pi_1\left(\Ms, +\infty\right)$.

Suppose now that $n$ is odd. Then
$f_0^{-1}\left(a^{-n}\right)[-\infty]=\ell_1^{-1}a^{(-n+1)/2}$. We
have
\[\widetilde{h_-}\left(f_0^{-1}\left(\gamma\right)[+\infty]\right)=
f_1^{-1}\left(\gamma\right)[-\infty]
\] for every $\gamma\in G_0$. Let us identify the group $\pi_1(\Ms,
+\infty)$ with $\pi_1\left(\Ms, -\infty\right)$, by identifying
the loop $f_1^{-1}\left(\gamma\right)[-\infty]\in\pi_1\left(\Ms,
-\infty\right)$ with the loop \[\ell_1\cdot
f_1^{-1}\left(\gamma\right)[-\infty]\cdot\ell_1^{-1}= \ell_1\cdot
f_0^{-1}\left(a^{-n}\gamma a^n\right)[-\infty]\cdot\ell_1^{-1}=
a^{-(n-1)/2}\cdot f_0^{-1}\left(\gamma\right)[+\infty]\cdot
a^{(n-1)/2}.\] We see that for every $\gamma\in\pi_1\left(\Ms,
+\infty\right)$ the loop $\widetilde{h_-}(\gamma)$ is identified
with $a^{-(n-1)/2}\gamma a^{(n-1)/2}$.

Therefore in all cases we can find a homeomorphism $h:\Ms\arr\Ms$
satisfying $f_0=f_1\circ h$ and such that the induced homomorphism
$h_*$ on the fundamental group $\pi_1\left(\Ms, +\infty\right)$ is
inner. But if a homeomorphism of a surface induces an inner
automorphism of the fundamental group, then this homeomorphism is
isotopic to identity. Consequently, $f_0$ and $f_1$ are homotopic.
\end{proof}

\subsection{The set $\F$ as a bimodule over the mapping class
group}\label{ss:bimodule}

Denote by $\mc$ the (pure) mapping class group of
$\Ms=\CS\setminus P$, i.e., the set of isotopy classes relative to
$P$ of homeomorphisms $h:\CS\arr\CS$ fixing $P$ pointwise. The
group $\mc$ is the quotient of the pure braid group on $|P|-1$
strings by its centre. For instance, if $|P|=4$, then $\mc$ is a
free group on two generators.

The set $\F$ has a natural structure of a permutational
$\mc$-bimodule, i.e.\ a set $\F$ equipped with commuting left- and
right-actions of the group $\mc$. If $f\in\F$ and $h\in\mc$, then we
just set $fh$ and $hf$ to be the corresponding compositions. (Here in
the composition $fh$ the map $f$ acts before $h$.) It is easy to see
that the left and the right actions of $\mc$ on $\F$ are well defined
and commute.

There exists a natural homomorphism
$\mc\arr\operatorname{Out}\left(\pi_1\left(\Ms\right)\right)$
mapping every element $h\in\mc$ to the automorphism $h_*$ of the
fundamental group of $\Ms$ (which is defined uniquely up to an
inner automorphism of $\pi_1(\Ms)$). It is well known
that this homomorphism is injective (see~\cite[Theorem~5.13.1]{zieschangVC}).

Using this fact and Proposition~\ref{pr:homotopy} one can describe
the structure of the bimodule $\F$. Take an arbitrary
$f\in\mathcal{F}$ and a homeomorphism $h:\CS\arr\CS$ acting
trivially on $P$ (we will also denote by $f\in\F$ and $h\in\mc$
the corresponding homotopy classes). 
Recall that $a$ denotes a small circle in the neighbourhood of
$\infty$.

\begin{proposition}
\label{pr:virten} For $f\in\F$ and $g, h\in\mc$, the following
conditions are equivalent:
\begin{itemize}
  \item[(i)]  $fg=hf$;
  \item[(ii)] there exists $n\in\Z$ such that
  \[\Phi_f\left(\gamma^{g^{-1}}\right)=\Phi_f\left(\gamma^{a^n}\right)^{h^{-1}}\]
for all $\gamma\in\pi_1\left(\Ms, +\infty\right)$;
  \item[(iii)] there exists $n\in\Z$ such that
   \[\Phi_f\left(\gamma\right)^h=\Phi_f\left(\gamma^{ga^n}\right)\]
  for all $\gamma\in\pi_1\left(\Ms, +\infty\right)$;
  \item[(iv)] there exists $n\in\Z$ such that
  \[\phi_f(\gamma)^h=\phi_f\left(
  \gamma^{ga^n}\right)\] for all $\gamma\in\dom\phi_f$, where $\phi_f$
  is the virtual endomorphism associated with the first coordinate of the wreath
  recursion $\Phi_f$.
\end{itemize}
\end{proposition}

\begin{proof}
  It is obvious that (ii) is equivalent to (iii). The equivalence of
  (iii) and (iv) follows directly from the definition of the virtual
  endomorphism associated with a wreath recursion. The equivalence of
  (i) and (ii) follows directly from
  Propositions~\ref{pr:homotopy},~\ref{pr:imgwr} and the specific
  definition of $\Phi_f$ which we gave in
  Subsection~\ref{ss:homotopy}.
\end{proof}

\begin{defi}
  Two branched coverings $f_1, f_2\in\mathcal{F}$ are
  \emph{combinatorially equivalent} (also said \emph{Thurston
    equivalent}) if there exists a homeomorphism $h$ of $\CS$ fixing
  $P$ pointwise, such that $f_1\cdot h$ and $h\cdot f_2$
  are homotopic.
\end{defi}

In other words, two elements $f_1, f_2\in\mathcal{F}$ are
equivalent if there exists an element $h$ of the mapping class
group of $\CS\setminus P$ such that $f_1h=hf_2$ in $\F$ (recall
that $\F$ is the set of \emph{homotopy classes} of branched
coverings).

\begin{corollary}
  \label{co:equivalence} If two branched coverings $f_0,
  f_1\in\mathcal{F}$ are combinatorially equivalent, then
  \[\Lambda_{f_0}\left(\pi_1\left(\Ms,
      t\right)\right)=\Lambda_{f_1}\left(\pi_1\left(\Ms,
      t\right)\right),\]
  as subsets of $\autxs$.
\end{corollary}
\noindent (Recall that $\Lambda_f:\pi_1(\Ms, t)\arr\autxs$ is the associated
iterated monodromy action.)
\begin{proof}
  Let first $\Phi:G\arr G\wr\symm$ be a wreath recursion, and let $h$
  be an automorphism of $G$. Consider the wreath recursion
  $\Phi^h:G\arr G\wr\symm$ given by $\Phi^h(g)=\Phi(g^{h^{-1}})^h$,
  where $h$ acts on $G\wr\symm$ by the diagonal action on $G^X$. Let
  $\Lambda_\Phi$ denote the action of $G$ on $\xs$ defined by the
  recursion $\Phi$: if $g\in G$, then $\Lambda_\Phi(g)$ maps $x_1\dots
  x_n\in\xs$ to $y_1\dots y_n$, where
  $\Phi(g)=\pair<g_1,\dots,g_d>\pi$ and $x_1^\pi=y_1$ and
  $\Lambda_\Phi(g_{x_1})$ maps $x_2\dots x_n$ to $y_2\dots y_n$.

  We easily check by induction that
  $\Lambda_\Phi(g)=\Lambda_{\Phi^h}(g^h)$: : the respective actions on
  the first level coincide, and the $x$th coordinate $g^h|_x$ of
  $\Phi^h(g^h)$ is equal to $(g|_x)^h$, where $g|_x$ is the $x$th
  coordinate of $\Phi(g)$.

  Let now $h$ be an element of the mapping class group of
  $\CS\setminus P_{f_0}$ such that $f_0\cdot h$ and $h\cdot f_1$ are
  homotopic. Then by Proposition~\ref{pr:homotopy}
  \[\Phi_{f_0}(\gamma^{h^{-1}})=\Phi_{f_1}(\gamma^{a^n})^{h^{-1}}\]
  for some $n$. Consequently,
  $\Phi_{f_0}^h(\gamma)=\Phi_{f_0}(\gamma^{h^{-1}})^h=\Phi_{f_1}(\gamma^{a^n})$.
  By the first two paragraphs, we have
  $\Lambda_{\Phi_{f_0}^h}(\gamma)=\Lambda_{f_0}(\gamma^{h^{-1}})$.

  The wreath recursion
  $\gamma\mapsto\Psi(\gamma):=\Phi_{f_1}(\gamma^{a^n})$ is the
  conjugate by $\Phi_{f_1}(a^n)$ of $\Phi_{f_1}$, so
  $\Lambda_\Psi(\gamma)$ is the conjugate of $\Lambda_{f_1}(\gamma)$
  by the automorphism $\Delta$ of $\xs$ given by the recursion
  \[\Delta=\tau^n\cdot\pair<\Delta, \Delta>,\]
  where $\tau=\pair<1,\tau>\sigma=\Lambda_{f_i}(a)$ is the adding machine. Now the element
  $\tau^n\Delta$ satisfies the recursion
  \[\tau^n\Delta)=\tau^{2n}\pair<\Delta,\Delta>=\pair<\tau^n\Delta,\tau^n\Delta>\]
  so $\tau^n\Delta=1$ and therefore $\Delta$ equals $\tau^{-n}$.

  Consequently,
  $\Lambda_{f_0}(\gamma^{h^{-1}})=\Lambda_{f_1}(\gamma^{a^{-n}})$,
  hence $\Lambda_{f_0}(\gamma)=\Lambda_{f_1}(\gamma^{ha^{-n}})$, which
  implies that $\Lambda_{f_0}(\pi_1(\Ms, t))=\Lambda_{f_1}(\pi_1(\Ms,
  t))$.
\end{proof}

\section{Twisted rabbits}\label{se:tr}

We consider now some concrete examples of bimodules $\F$. We will
consider the cases when $f$ is a quadratic polynomial whose set of
finite post-critical points has size $3$.

Let us consider first the ``Douady
rabbit''~\cite[Figure~35]{milnor}. It is the polynomial
$f_R\approx z^2+(-0.1226 + 0.7449i)$. The two other polynomials
inducing the same permutation of their post-critical set are the
``corabbit'' $f_C\approx z^2+(-0.1226 - 0.7449i)$ and the
``airplane'' $f_A\approx z^2-1.7549$.

\begin{figure}[ht]
\centering
  \includegraphics{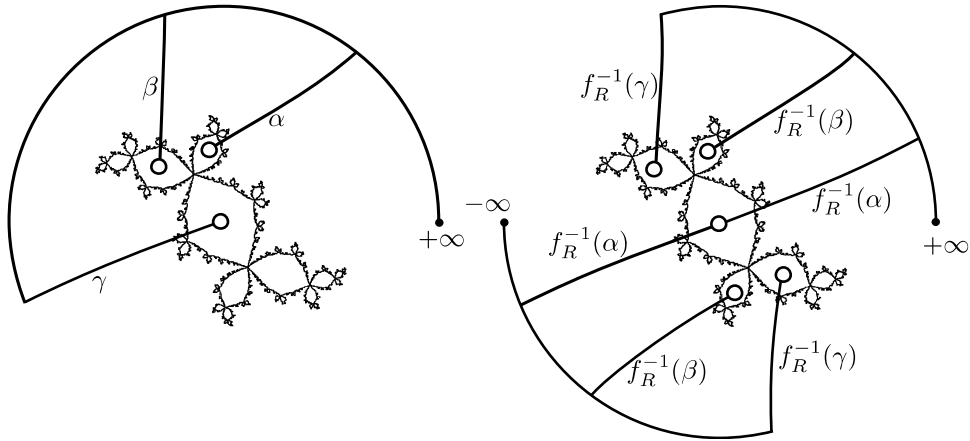}
  \caption{Computing $\img{f_R}$}\label{fig:rabbitl}
\end{figure}

We choose as usual $+\infty$ as the basepoint. Let $\alpha, \beta,
\gamma$ be the loops going around, respectively, $c, c^2+c$ and
$0$ in the positive direction and connected to the basepoint as it
is shown on left-hand side part of Figure~\ref{fig:rabbitl}. Let
$P=\{\infty, 0, c, c^2+c\}$ be the post-critical set of $f_R$.

\noindent The ``rabbit'''s wreath recursion $\Phi_{f_R}$ is defined by
\[
\Phi_{f_R}(\alpha)=
\pair<\alpha^{-1}\beta^{-1},\gamma\beta\alpha>\sigma,\quad
\Phi_{f_R}(\beta)=\pair<\alpha,1>,\quad\Phi_{f_R}(\gamma)=\pair<\beta,1>.
\]
The preimages of the paths $\alpha, \beta, \gamma$ are shown on the
right-hand side part of Figure~\ref{fig:rabbitl}.  Note that
$\tau=\gamma\beta\alpha=\pair<1,\gamma\beta\alpha>\sigma$ is the
standard adding machine.

\subsection{The mapping class group action}
\begin{figure}[ht]
\centering\includegraphics{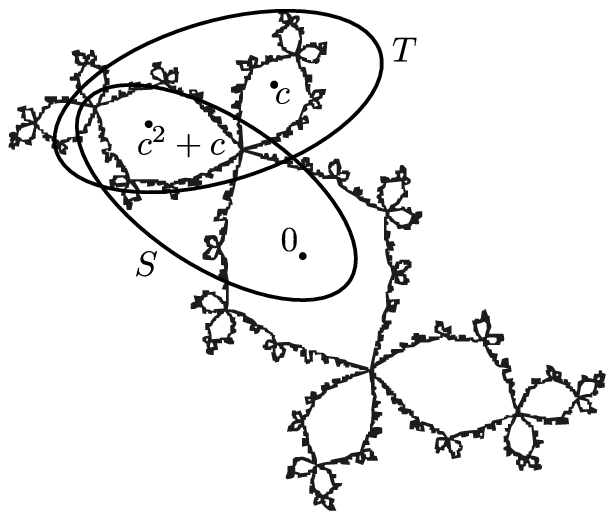} \caption{The generators of
$\mc$}\label{fig:st}
\end{figure}

The mapping class group $\mc$ is generated by the left-handed
(counter-clockwise) Dehn twist $T$ about the curve encircling the
points $c, c^2+c$ and by the left-handed twist $S$ about the curve
encircling the points $0, c^2+c$ (see Figure~\ref{fig:st}).

The twists $T$ and $S$ are defined by their action
on the fundamental group of the punctured plane by the rules
\begin{xalignat*}{3}
  \alpha^T&=\alpha^{\beta\alpha},&\beta^T&=\beta^{\beta\alpha}=
  \beta^{\alpha},&\gamma^T&=\gamma,\\
  \alpha^S&=\alpha,&\beta^S&=\beta^{\gamma\beta},&
  \gamma^S&=\gamma^{\gamma\beta}=\gamma^{\beta}.\end{xalignat*}
Their inverses act by the rules
\begin{xalignat*}{3}
\alpha^{T^{-1}}&=\alpha^{\alpha^{-1}\beta^{-1}}=\alpha^{\beta^{-1}},&
\beta^{T^{-1}}&=\beta^{\alpha^{-1}\beta^{-1}},&
\gamma^{T^{-1}}&=\gamma,\\
\alpha^{S^{-1}}&=\alpha,&
\beta^{S^{-1}}&=\beta^{\beta^{-1}\gamma^{-1}}=\beta^{\gamma^{-1}},&
\gamma^{S^{-1}}&=\gamma^{\beta^{-1}\gamma^{-1}}.\end{xalignat*}

\begin{proposition}
\label{pr:psibar} Let $\psi$ be the virtual endomorphism of the
group $\mc$ defined on the subgroup $H=\left\langle T^2, S,
S^T\right\rangle$ of index $2$ by
\[\psi\left(T^2\right)=S^{-1}T^{-1},\quad \psi\left(S\right)=T,\quad
\psi\left(S^T\right)=1.\] Consider the map
\[\opsi:g\mapsto\begin{cases}\psi(g) & \text{ if $g$ belongs to the
domain of $\psi$},\\ T\psi(gT^{-1}) & \text{
otherwise.}\end{cases}\] Then for every $g\in\mc$ the branched
coverings $f_R\cdot g$ and $f_R\cdot\opsi(g)$ are combinatorially
equivalent.
\end{proposition}

The subgroup $H$ is generated by those Dehn twists about curves that
encircle an even number of times the critical value $c$. They are
therefore those mapping classes that can be lifted through $f_R$. The
map $\psi$ is precisely that lift.

\begin{proof} We first claim that for all $g\in H$ the maps $f_Rg$ and
  $\psi(g)f_R$ are homotopic. It suffices to check this on the
  generators $\{T^2, S, S^T\}$ of $H$, and this is done below.

  Note then that $f_R\cdot g=\psi(g)\cdot
  f_R=\psi(g)\cdot\left(f_R\cdot\psi(g)\right)\cdot\psi(g)^{-1}$,
  i.e., that $f_R\cdot g$ and $f_R\cdot\psi(g)$ are combinatorially
  equivalent.

  If $g$ does not belong to the domain of $\psi$, then $gT^{-1}$ does,
  and \[f_R\cdot g=f_R\cdot gT^{-1}T=\psi\left(gT^{-1}\right)\cdot
  f_R\cdot T= \psi\left(gT^{-1}\right)\cdot\left(f_R\cdot
    T\psi\left(gT^{-1}\right)\right)\cdot\psi\left(gT^{-1}\right)^{-1},\]
  i.e., $f_R\cdot g$ and $f_R\cdot\opsi(g)$ are combinatorially
  equivalent.

  Let now $\phi=\phi_{f_R}$ be the virtual endomorphism associated with
  the first coordinate of the wreath recursion $\Phi_{f_R}$. We have
\begin{gather*}
\alpha =\phi\left(\beta\right),\quad \beta
=\phi\left(\gamma\right),\quad \gamma
=\phi\left(\alpha^{2\beta^{-1}\gamma^{-1}}\right),\\
\gamma^{\beta\alpha} = \phi\left(\alpha^2\right),\quad
\phi\left(\beta^\alpha\right)=1,\quad\phi\left(\gamma^\alpha\right)=1,\\
\intertext{from which we compute}
\phi\left(\beta^{T^2}\right)=\phi\left(\beta^{\alpha\beta\alpha}\right)=
\phi\left(\beta^{\alpha^2\cdot\beta^\alpha}\right)=\alpha^{\beta^{-1}\gamma\beta\alpha},\\
\phi\left(\gamma^{T^2}\right)=\phi\left(\gamma\right)=\beta,\\
\phi\left(\left(\alpha^{2\beta^{-1}\gamma^{-1}}\right)^{T^2}\right)=
\phi\left(\alpha^{2\beta^{-1}\beta\alpha\beta\alpha\gamma^{-1}}\right)=
\phi\left(\alpha^{2\beta^\alpha\cdot\gamma^{-1}}\right)=
\gamma^{\beta\alpha\beta^{-1}}.
\end{gather*}
We see that for every $\delta\in\pi_1\left(\Ms, +\infty\right)$ we
have $\phi\left(\delta^{T^2}\right)=\phi(\delta)^h$, where the
automorphism $h$ is given on the generators by
\begin{align}
\alpha^h &= \alpha^{\beta^{-1}\gamma\beta\alpha},\notag\\
\label{eq:h123} \beta^h  &= \beta,\\
\gamma^h &= \gamma^{\beta\alpha\beta^{-1}}.\notag
\end{align}
The automorphism $h$ is equal to the product $S^{-1}T^{-1}a$,
where $a$ is conjugation by $\gamma\beta\alpha$:
\[\begin{array}{cclclcl}
\alpha &\stackrel{S^{-1}}{\mapsto} & \alpha &
\stackrel{T^{-1}}{\mapsto} & \alpha^{\beta^{-1}} &
\stackrel{a}{\mapsto} &
\alpha^{\beta^{-1}\gamma\beta\alpha},\\
\beta &\stackrel{S^{-1}}{\mapsto} & \beta^{\gamma^{-1}}
&\stackrel{T^{-1}}{\mapsto}
& \beta^{\alpha^{-1}\beta^{-1}\gamma^{-1}} &\stackrel{a}{\mapsto} & \beta,\\
\gamma &\stackrel{S^{-1}}{\mapsto} &
\gamma^{\beta^{-1}\gamma^{-1}} &\stackrel{T^{-1}}{\mapsto} &
\gamma^{\beta^{-\alpha^{-1}\beta^{-1}}\gamma^{-1}}
&\stackrel{a}{\mapsto} & \gamma^{\beta\alpha\beta^{-1}}.
\end{array}\]
Consequently, by Proposition~\ref{pr:virten} (iv) with $n=1$,
\[f_R\cdot T^2=S^{-1}T^{-1}\cdot f_R.\]

We have next
\begin{gather*}
\phi\left(\beta^S\right)=\phi\left(\beta^{\gamma\beta}\right)=\alpha^{\beta\alpha},\\
\phi\left(\gamma^S\right)=\phi\left(\gamma^\beta\right)=\beta^\alpha,\\
\phi\left(\left(\alpha^{2\beta^{-1}\gamma^{-1}}\right)^S\right)=
\phi\left(\alpha^{2\beta^{-1}\gamma^{-1}}\right)=\gamma,\\
\intertext{so} f_R\cdot S=T\cdot f_R.
\end{gather*}

Let us finally compute the action of $S^T$ on the generators
$\{\alpha, \beta, \gamma\}$ of the fundamental group:
\[\begin{array}{lclclcl}
\alpha &\stackrel{T^{-1}}{\mapsto} & \alpha^{\beta^{-1}} &
\stackrel{S}{\mapsto} & \alpha^{\beta^{-\gamma\beta}} &
\stackrel{T}{\mapsto} &
\alpha^{\beta\alpha\beta^{-\alpha\gamma\beta^\alpha}}=\alpha^{\gamma^{-1}\alpha^{-1}\beta^{-1}\alpha\gamma\alpha^{-1}\beta\alpha},\\
 \beta & \stackrel{T^{-1}}{\mapsto} &
\beta^{\alpha^{-1}\beta^{-1}} & \stackrel{S}{\mapsto} &
\beta^{\gamma\beta\alpha^{-1}\beta^{-1}\gamma^{-1}\beta^{-1}\gamma\beta}
& \stackrel{T}{\mapsto} &
\beta^{\alpha\gamma\alpha^{-1}\gamma^{-1}\alpha^{-1}
\beta^{-1}\alpha\gamma\alpha^{-1}\beta\alpha},\\
\gamma & \stackrel{T^{-1}}{\mapsto} & \gamma &
\stackrel{S}{\mapsto} & \gamma^\beta & \stackrel{T}{\mapsto} &
\gamma^{\alpha^{-1}\beta\alpha}.
\end{array}\]

Let us apply the virtual endomorphism $\phi$ to the conjugators:
\begin{gather*}
\phi\left(\gamma^{-1}\cdot\alpha^{-1}\beta^{-1}\alpha\cdot\gamma
\cdot\alpha^{-1}\beta\alpha\right)=\beta^{-1}\beta=1,\\
\phi\left(\alpha\gamma\alpha^{-1}\cdot\gamma^{-1}\cdot\alpha^{-1}
\beta^{-1}\alpha\cdot \gamma\cdot\alpha^{-1}\beta\alpha\right)=
\beta^{-1}\beta=1,\\
\phi\left(\alpha^{-1}\beta\alpha\right)=1.\\
\intertext{Consequently,}
\phi\left(\left(\alpha^2\right)^{T^{-1}ST}\right)=\phi\left(\alpha^2\right),\quad
\phi\left(\beta^{T^{-1}ST}\right)=\phi\left(\beta\right),\quad
\phi\left(\gamma^{T^{-1}ST}\right)=\phi\left(\gamma\right)\\
\intertext{which implies that}f_R\cdot S^T=f_R.\qedhere
\end{gather*}
\end{proof}

\begin{proposition}
\label{pr:opsicontracting} The map $\opsi:\mc\arr\mc$ is
contracting: for every $g\in\mc$ there exists $n\in\mathbb{N}$
such that $\opsi^{n}(g)\in\{1, T, T^{-1}\}$.
\end{proposition}

\begin{proof}
  Consider the wreath recursion for $\mc$ given by
  \begin{equation}\label{eq:mcrec}
  \Phi(T)=\pair<1,S^{-1}T^{-1}>\sigma,\quad\Phi(S)=\pair<T,1>.\end{equation}
  It is straightforward to check that $\opsi$ may be recovered from
  the recursion as follows: if $\Phi(g)=\pair<g_0,g_1>$, then
  $\opsi(g)=g_0$; if $\Phi(g)=\pair<g_0,g_1>\sigma$, then
  $\opsi(g)=Tg_0$. Since $T|_0=1$, we obtain inductively for all
  $n\in\N$
  \[\opsi^n(g)=g|_v\text{ or }Tg|_v\text{ for some }v\in X^n.\]

  We first claim that the recursion~\eqref{eq:mcrec} is
  contracting. For that purpose, it suffices to compute the nucleus
  of~\eqref{eq:mcrec}. A simple calculation shows that it is
  \[\nuke=\{1,S,T,TS,S^{-1},T^{-1},S^{-1}T^{-1}\}.\]
  It follows that for all $g\in\mc$ we have $\opsi^n\in\nuke\cup
  T\nuke$ for all sufficiently large $n$, and therefore that $\opsi^n(g)$
  lands on a $\opsi$-cycle. Direct computations then show that the only
  $\opsi$-periodic elements in $\nuke\cup T\nuke$ are the fixed points $1$, $T$
  and the cycle $T^{-1}\arr T^2S\arr S^{-1}\arr T^{-1}$.
\end{proof}

\subsection[Contraction along the subgroup $\langle T\rangle$]{Contraction along the subgroup $\boldsymbol{\langle T\rangle}$}
Every integer $m$ has a unique \emph{4-adic} expansion
\[m=\sum_{k=0}^\infty m_k4^k,\] with $m_k\in\{0,1,2,3\}$, and almost
all $m_k=0$ if $m$ is non-negative, and almost all $m_k=3$ if $m$ is
negative.\footnote{For example, $m=-1$ corresponds to $m_k=3$ for all
  $k$.}

\begin{proposition}
\label{pr:4adic} If the 4-adic expansion of the number $m$ has
digits $1$ or $2$, then the branched covering $f_R\cdot T^m$ is
equivalent to the $f_R\cdot T$. Otherwise it is equivalent to
$f_R$ for non-negative $m$ and to $f_R\cdot T^{-1}$ for negative
$m$.
\end{proposition}

\begin{proof}
Let us iterate the map $\opsi$ on the cyclic subgroup
$\left\langle T\right\rangle$. We have
\begin{align*}
\opsi^3(T^{4k})&=\opsi^2\left(S^{-1}T^{-1}\right)^{2k}
=\opsi^2\left(\left(S^{-1}\cdot S^{-T}\cdot T^{-2}
\right)^k\right)\\&=\opsi\left(T^{-1}\cdot 1\cdot
TS\right)^k=\opsi\left(S^k\right)=T^k,
\end{align*} so $f_R\cdot T^{4k}$ is equivalent to
$f_R\cdot T^k$. Similarly
\begin{align*}
\opsi^3(T^{4k+1})&=\opsi^2\left(T\psi\left(T^{4k}\right)\right)
=\opsi^2\left(T\left(S^{-1}T^{-1}\right)^{2k}\right)\\
&=\opsi\left(T\psi\left(T\left(S^{-1}T^{-1}\right)^{2k}T^{-1}\right)\right)=
\opsi\left(T\psi\left(\left(TS^{-1}T^{-1}\cdot S^{-1}\cdot
T^{-2}\right)^k\right)\right)\\
&=\opsi\left(T\left(1\cdot T^{-1}\cdot
TS\right)^k\right)=\opsi\left(TS^k\right)\\
&=T\psi\left(TS^kT^{-1}\right)=T,
\end{align*} so all branched coverings $f_R\cdot T^{4k+1}$ are equivalent to
$f_R\cdot T$. Next
\begin{align*}
\opsi^3(T^{4k+2})&=\opsi^2\left(S^{-1}T^{-1}\right)^{2k+1}
=\opsi\left(T\psi\left(\left(S^{-1}T^{-1}\right)^{2k+1}T^{-1}\right)\right)\\
&=\opsi\left(T\psi\left(\left(S^{-1}\cdot S^{-T}\cdot T^{-2}\right)^k\cdot S^{-1}\cdot T^{-2}\right)\right)\\
&=\opsi\left(T \left(T^{-1}\cdot 1\cdot TS\right)^k\cdot
T^{-1}\cdot TS\right)=\opsi\left(TS^{k+1}\right)\\
&=T\psi\left(TS^{k+1}T^{-1}\right)=T,\end{align*} so all branched
coverings $f_R\cdot T^{4k+2}$ are equivalent to $f_R\cdot T$.
Finally
\begin{align*}
\opsi^3(T^{4k+3})&=\opsi^2\left(T\psi\left(T^{4k+2}\right)\right)
=\opsi^2\left(T\left(S^{-1}T^{-1}\right)^{2k+1}\right)\\
&=\opsi^2\left(TS^{-1}T^{-1}\cdot\left(S^{-1}\cdot S^{-T}\cdot
T^{-2}\right)^k\right)=\opsi\left(1\cdot\left(T^{-1}\cdot 1\cdot
TS\right)^k\right)\\ &=\opsi\left(S^k\right)=T^k,\end{align*} so
$f_R\cdot T^{4k+3}$ is equivalent to $f_R\cdot T^k$. The statement
of the proposition now easily follows.
\end{proof}

\subsection[Solving the problem for all $m\in\Z$]{Solving the problem for all $\boldsymbol{m\in\Z}$}\label{se:m}

\begin{figure}[ht]
{\centering
  \includegraphics{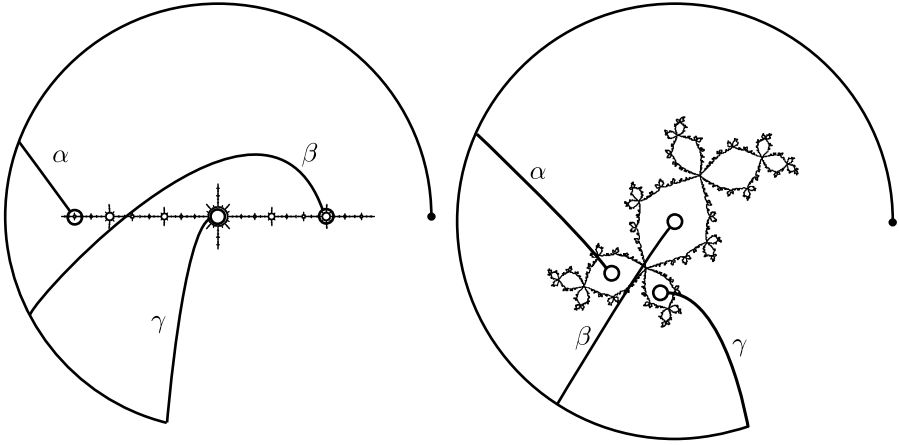}
}\caption{The standard generators of $\img{f_A}$ and
$\img{f_C}$}\label{fig:caray}
\end{figure}

The wreath recursions for the ``airplane'' polynomial $f_A$ is
given by
\begin{equation}\label{eq:recairplane}
\Phi_{f_A}\left(\alpha\right)=\pair<\alpha^{-1},\gamma\alpha>\sigma,\quad
\Phi_{f_A}\left(\beta\right)=\pair<\alpha,1>,\quad
\Phi_{f_A}\left(\gamma\right)=\pair<1,\beta^{\gamma^{-1}}>.\end{equation}
Here $\alpha, \beta, \gamma$ are loops going in the positive
direction around $c, c^2+c, (c^2+c)^2+c$ respectively, as shown on
the left part of Figure~\ref{fig:caray}.

The wreath recursion for the ``corabbit'' $f_C$ is given by
\begin{equation}\label{eq:reccorabbit}
\Phi_{f_C}\left(\alpha\right)=\pair<\alpha^{-1}\beta^{-1},\gamma\beta\alpha>\sigma,
\quad
\Phi_{f_C}\left(\beta\right)=\pair<\alpha^{\beta\alpha},1>,\quad
\Phi_{f_C}\left(\gamma\right)=\pair<\beta^\alpha,1>.\end{equation}
Here $\alpha, \beta, \gamma$ are loops going in the positive
direction around $c, c^2+c, (c^2+c)^2+c$ respectively, as shown on
the right part of Figure~\ref{fig:caray}. Note that for $f=f_A,
f_C$, just as for $f=f_R$, we have
$\Phi_f(\gamma\beta\alpha)=\pair<1, \gamma\beta\alpha>\sigma$.

Let us identify the planes of $f_A$ and $f_C$ with the plane of
$f_R$ by identifying their respective loops $\alpha, \beta,
\gamma$ (the definition of $\alpha, \beta, \gamma$ in the plane of
$f_R$ is given on the left-hand side of Figure~\ref{fig:rabbitl}).

Let $T$ denote, as before, the left Dehn twist about the curve around
the points $c, c^2+c$ in the plane of $f_R$. Then, from $T$'s
definition, we get the following wreath recursion for the standard
iterated monodromy action for $T^m\cdot f_R$:
\begin{equation}\label{eq:twistedrec}
\begin{array}{l}\Phi_{T^m\cdot
f_R}\left(\alpha\right)=\Phi_{f_R}\left(\alpha\right)^{T^{-m}}
=\pair<\alpha^{-1}\beta^{-1},\gamma\beta\alpha>\sigma,\\
\Phi_{T^m\cdot
f_R}\left(\beta\right)=\Phi_{f_R}\left(\beta\right)^{T^{-m}}=
\pair<\alpha^{\left(\alpha^{-1}\beta^{-1}\right)^m},1>,\\
\Phi_{T^m\cdot
f_R}\left(\gamma\right)=\Phi_{f_R}\left(\gamma\right)^{T^{-m}}=
\pair<\beta^{\left(\alpha^{-1}\beta^{-1}\right)^m},1>,
\end{array}\end{equation}
where $T^{-m}$ acts on $\pi_1\left(\Ms, +\infty\right)\wr\symm$ by
the diagonal action:
\[\left(\pair<x,y>\sigma^k\right)^{T^{-m}}=\pair<x^{T^{-m}},y^{T^{-m}}>\sigma^k.\]

Corollary~\ref{co:equivalence} makes it possible to solve Hubbard's
question algorithmically for every given $m$ in the following way.

Thurston's Theorem~\ref{th:thurston} implies that $T^m\cdot f_R$
is combinatorially equivalent to exactly one polynomial in the set
$\{f_R, f_A, f_C\}$. There are no obstructions, since the only
obstructions for polynomials are Levy cycles, which cannot exist
in the case of a periodic critical point.
Corollary~\ref{co:equivalence} then tells us that
$\Lambda_{T^m\cdot
  f_R}\left(\pi_1\left(\Ms\right)\right)$ coincides with the iterated
monodromy group of the associated polynomial. One can prove that these
groups are different (as sets), and therefore, if we prove that
$\Lambda_{T^m\cdot f_R}\left(\pi_1\left(\Ms\right)\right)$ coincides
with a given group $\img{f_*}$ for $*\in\{R,A,C\}$, then we can
conclude that $T^m\cdot f_R$ is equivalent to the respective $f_*$.

\begin{figure}[ht]
\centering
  \includegraphics{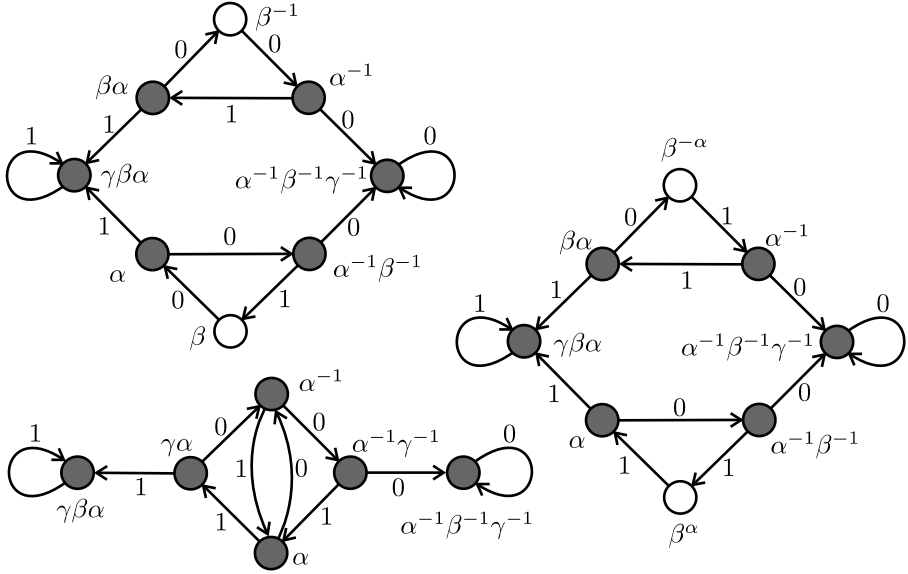}
  \caption{Nuclei of the ``rabbit'' (top), the ``corabbit'' (right) and
  the ``airplane'' (bottom)}\label{fig:rabnuc}
\end{figure}

We therefore prove that the $\img{f_*}$ are all distinct. This is
done by computing their nuclei, and checking that they are
distinct as finite automata; this is done in
Figure~\ref{fig:rabnuc}.

\begin{proposition}\label{prop:R1=airplane}
The group $\img{T\cdot f_R}=\Lambda_{T\cdot
f_R}\left(\pi_1\left(\Ms\right)\right)$ coincides with
$\img{f_A}$. Indeed the homeomorphism $h=TS^{-1}a$ conjugates
$T\cdot f_R$ with $f_A$, if the planes of $f_R,f_A$ are identified
as above.
\end{proposition}

\begin{proof} Let $\alpha, \beta, \gamma$ be the generators of
$\img{f_A}$. They are defined now as the automorphisms of $\xs$
satisfying the recursion (compare with~\eqref{eq:recairplane})
\[\alpha=\pair<\alpha^{-1},\gamma\alpha>\sigma,\quad
\beta=\pair<\alpha,1>,\quad \gamma=\pair<1,\beta^{\gamma^{-1}}>.\]

Let $\alpha_1, \beta_1$ and $\gamma_1$ be the generators of
$\img{T\cdot f_R}$. They are given by the
recursion~\eqref{eq:twistedrec}:
\[\alpha_1=\pair<\alpha_1^{-1}\beta_1^{-1},\gamma_1\beta_1\alpha_1>\sigma,\quad
\beta_1=\pair<\alpha_1^{\beta_1^{-1}},1>,\quad
\gamma_1=\pair<\beta_1^{\alpha_1^{-1}\beta_1^{-1}},1>.\]

We claim that
\[\alpha_1=\alpha^h=\alpha^{\beta^{\gamma^{-1}}\alpha\gamma\beta\alpha},\quad
\beta_1=\beta^h=\beta^{\gamma^{-1}\alpha\gamma\beta\alpha},\quad
\gamma_1=\gamma^h=\gamma^\alpha.\] For that purpose, it suffices
to show that the left-hand sides of these equalities satisfy the
same recursions as $\alpha_1, \beta_1, \gamma_1$. Note that
$\gamma^h\beta^h\gamma^h=(\gamma\beta\alpha)^h=\gamma\beta\alpha$.
We have
\begin{align*}
\beta^{\gamma^{-1}}\alpha\gamma\beta\alpha&=\pair<\alpha,1>\pair<\alpha^{-1},\gamma\alpha>\sigma
\pair<1, \gamma\beta\alpha>\sigma=\pair<\gamma\beta\alpha,
\gamma\alpha>,\\
\intertext{hence}
\alpha^{\beta^{\gamma^{-1}}\alpha\gamma\beta\alpha}&=
\left(\pair<\alpha^{-1},\gamma\alpha>\sigma\right)^{\pair<\gamma\beta\alpha,\gamma\alpha>}
=\pair<\alpha^{-1}\beta^{-1}\gamma^{-1}\alpha^{-1}\gamma\alpha,
\alpha^{-1}\gamma^{-1}\gamma\alpha\gamma\beta\alpha>\sigma\\
&=\pair<\alpha^{-1}\beta^{-1}\gamma^{-1}\alpha^{-1}\gamma\alpha,\gamma\beta\alpha>\sigma
=\pair<(\gamma\beta\alpha)^{-h}\gamma^h,
(\gamma\beta\alpha)^h>\sigma\\
&=\pair<\alpha^{-h}\beta^{-h},\gamma^h\beta^h\alpha^h>\sigma,\\
\beta^{\gamma^{-1}\alpha\gamma\beta\alpha}&=\beta^{\gamma^{-1}\beta^{\gamma^{-1}}\alpha\gamma\beta\alpha}=
\left(\pair<\alpha,1>\right)^{\pair<\gamma\beta\alpha,\gamma\alpha>}=
\pair<\alpha^{\gamma\beta\alpha},1>\\
&=\pair<\left(\alpha^{\gamma\beta\gamma^{-1}\alpha\gamma\beta\alpha}\right)
^{\beta^{-\gamma^{-1}\alpha\gamma\beta\alpha}},1>=\pair<\left(\alpha^h\right)^{\beta^{-h}},1>,\\
\gamma^\alpha&=\left(\pair<1,\beta^{\gamma^{-1}}>\right)^{\pair<\alpha^{-1},\gamma\alpha>\sigma}
=\pair<\beta^{\alpha},1>=\pair<\left(\beta^{\alpha^{-\gamma\beta}\beta^{-1}}\right)^{\gamma^{-1}\alpha\gamma\beta\alpha},1>\\
&=\pair<\left(\beta^h\right)^{\alpha^{-h}\beta^{-h}},1>.\end{align*}

This shows the required relations between $\alpha, \beta, \gamma$
and $\alpha_1, \beta_1, \gamma_1$, so $\img{f_A}=\img{T\cdot
f_R}$. Proposition~\ref{pr:homotopy} now implies that $T\cdot f_R$
is homotopic to $h^{-1}\cdot f_A\cdot h$.
\end{proof}

\begin{proposition}
The group $\img{T^{-1}\cdot f_R}=\Lambda_{T^{-1}\cdot
f_R}\left(\pi_1\left(\Ms\right)\right)$ coincides with
$\img{f_C}$. Moreover, $T^{-1}\cdot f_R$ and $f_C$ are homotopic
if the planes of $f_R,f_C$ are identified as above.
\end{proposition}

\begin{proof}
Let $\alpha, \beta, \gamma$ be the generators of the iterated
monodromy group of the ``corabbit''. They are defined by the
recursion (compare with~\eqref{eq:reccorabbit})
\[\alpha=\pair<\alpha^{-1}\beta^{-1},\gamma\beta\alpha>\sigma,
\quad \beta=\pair<\alpha^{\beta\alpha},1>,\quad
\gamma=\pair<\beta^\alpha,1>.\]

Let $\alpha_{-1}, \beta_{-1}$ and $\gamma_{-1}$ be the generators
of $\img{T^{-1}\cdot f_R}$. They are given by the
recursion~\eqref{eq:twistedrec}:
\[\alpha_{-1}=\pair<\alpha_{-1}^{-1}\beta_{-1}^{-1},\gamma_{-1}\beta_{-1}\alpha_{-1}>\sigma,\quad
\beta_{-1}=\pair<\alpha_{-1}^{\beta_{-1}\alpha_{-1}},1>,\quad\gamma_{-1}
=\pair<\beta_{-1}^{\alpha_{-1}},1>.\] Since these two recursions
are the same, we have $\alpha_{-1}=\alpha$, $\beta_{-1}=\beta$,
$\gamma_{-1}=\gamma$. Proposition~\ref{pr:homotopy} now implies
that $f_C$ and $T^{-1}\cdot f_R$ are homotopic.
\end{proof}

\begin{corollary}
The branched coverings $T\cdot f_R$ and $f_R\cdot T$ are
equivalent to $f_A$ and the branched coverings $T^{-1}\cdot f_R$
and $f_R\cdot T^{-1}$ are equivalent to $f_C$.
\end{corollary}

The last corollary together with Proposition~\ref{pr:4adic} prove
the following solution of the ``twisted rabbit question''.

\begin{theorem}
\label{th:4adic} If the 4-adic expansion of the number $m$ has
digits $1$ or $2$, then the branched covering $f_R\cdot T^m$ is
equivalent to the ``airplane'' $f_A$. Otherwise it is equivalent
to the ``rabbit'' $f_R$ for non-negative $m$ and to the
``corabbit'' $f_C$ for negative $m$.
\end{theorem}

The general case of any element of the mapping class group $\mc$ is
treated using Proposition~\ref{pr:opsicontracting} in the following
theorem.

\begin{theorem}
\label{th:generalrabbit}
Let $g\in\mc$ be an arbitrary homeomorphism fixing the post-critical
set pointwise. Let $\opsi:\mc\arr\mc$ be the map defined in
Proposition~\ref{pr:psibar}. The orbit of $g$ under the iteration of $\opsi$
will land either on $1$, or on $T$, or on $T^{-1}$. In the first
case $g\cdot f_R$ is equivalent to the ``rabbit'', in the second case
it is equivalent to the ``airplane'' and in the last to the ``corabbit''.
\end{theorem}

This theorem gives an algorithm solving the general ``twisted rabbit
question''. Note that due to the fact that $\opsi$ is contracting on
$\mc$, this algorithm has linear complexity with respect to the
word-length of elements of $\mc$.

Note also that in Theorem~\ref{th:4adic} the typical answer is
``airplane'', and exponentially few values of $m$ yield ``rabbits'' or
``corabbits''. This seems to happen quite often on cyclic subgroups of
$\mc$. On the other hand, on the cyclic subgroup $\langle ST^2\rangle$
all twists are ``rabbits''; and on the subgroup $\langle ST\rangle$
there are roughly as many ``rabbits'' as ``airplanes'': $(ST)^m\cdot
f_R$ is equivalent to

\begin{center}
\begin{tabular}{ll}
  $f_R$ & if $m<-2,m=2n-3,$ and $n$'s 4-adic expansion has
    only 0's and 3's,\\
  $f_C$ & if $m>-2,m=2n-1,$ and $n$'s 4-adic expansion has
    only 0's and 3's,\\
  $f_R$ & if $n\equiv 0\pmod2$,\\
  $f_A$ & if $n\equiv 1\pmod2$ in the cases not covered above.
\end{tabular}
\end{center}

This is because $\opsi((ST)^{2m})=S^{-mT^{-1}}$ and
$\opsi((ST)^{2m+1})=T^2S^{-m}$; calculations are similar to those of
Proposition~\ref{pr:4adic}.

These phenomena can ultimately be traced to the following reason: the
subgroup $\dom\psi^3$ of $\mc$ has index $8$:
\[\dom\psi^3=\langle S^4,T^4,ST^2,S^{-2}T^2S,S^{-1}T^2S^2,S^T,S^{T^{-1}},S^{TS},S^{T^{-1}S}\rangle.\]
The map $\psi^3$ is defined on these generators by $\psi^3(S^4)=S$,
$\psi^3(T^4)=T$, and all other generators are mapped to $1$.

\section{Dynamics on moduli space}\label{se:moduli}

We have computed in the previous sections a virtual endomorphism
$\psi$ of the mapping class group $\mc$ associated with the bimodule
$\F$ of branched coverings with the ramification graph of the
``rabbit'' polynomial. This virtual endomorphism is given by
\[\psi(T^2)=S^{-1}T^{-1},\quad\psi(S)=T,\quad\phi(S^T)=1.\]

This virtual endomorphism is associated to a self-similar action
of $\mc$ on the binary tree. Let us show that this action
coincides with the iterated monodromy action associated with the
rational function $F(w)=1-\frac{1}{w^2}$. 


The critical points of $1-\frac{1}{w^2}$ are $0$ and $\infty$.
Their orbit is
\[0\mapsto\infty\mapsto 1\mapsto 0.\]

\begin{figure}[ht]
\centering
  \includegraphics{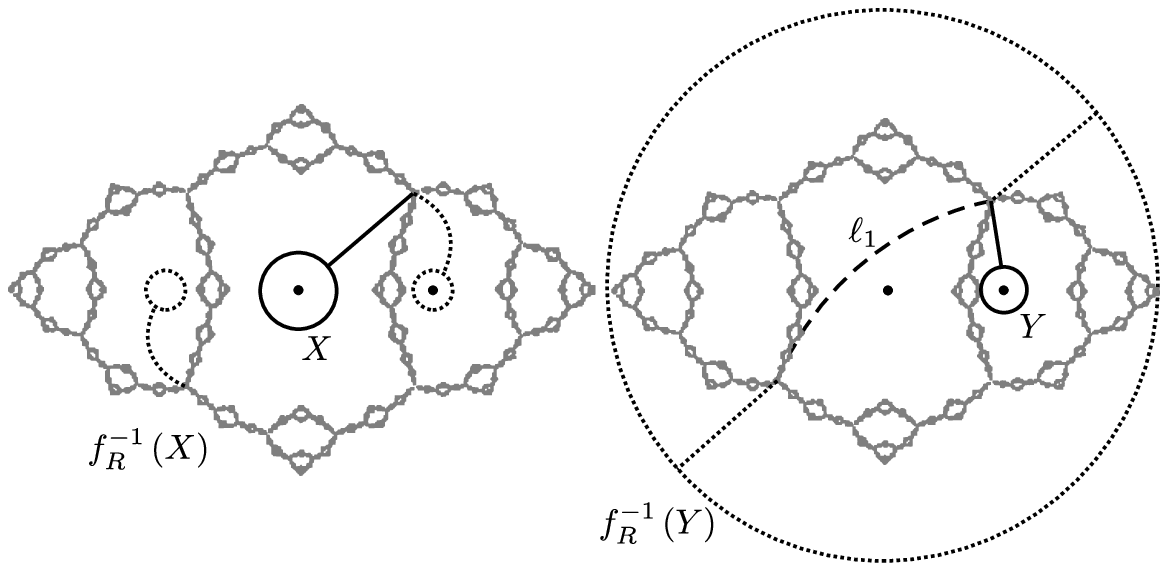}  \caption{}\label{fig:x}
\end{figure}

Let us take the fixed point $t\approx 0.8774 + 0.7449i$ as our
basepoint. The fundamental group of $\CS\setminus\{0, 1, \infty\}$
is generated by the loops going in the \emph{negative} direction
around $0$ and $1$. Let us denote the first by $X$ and the second
by $Y$ (see the solid lines on the left- and right-hand sides of
Figure~\ref{fig:x} respectively). Let us choose the connecting
paths $\ell_0$ and $\ell_1$, so that $\ell_0$ is the trivial path
at $t$ and $\ell_1$ connects $t$ to $-t$ passing above the
puncture $0$, as shown in the right-hand side of Figure~\ref{fig:x}.

We then get the following recursion for the iterated monodromy
group of $F$:
\[X=\pair<Y,1>,\quad Y =\pair<1,X^{-1}Y^{-1}>\sigma.\]
Next $Y^2=\pair<X^{-1}Y^{-1},X^{-1}Y^{-1}>$ and $Y^{-1}XY=
\pair<1,Y>$. We see that the virtual endomorphism associated with
the wreath recursion is given by
\[X\mapsto Y,\quad Y^2\mapsto X^{-1}Y^{-1},\quad Y^{-1}XY\mapsto 1.\]

It is therefore precisely the virtual endomorphism $\psi$
associated to $\F$, if we identify $X$ with $S$ and $Y$ with $T$.

\subsection{Moduli and Teichm\"uller space}
This coincidence has a nice explanation in terms of Teichm\"uller
theory. Let $P\subset \Sph$ be a finite subset of the sphere. Then
the \emph{Teichm\"uller space} $\mathcal{T}_P$ modelled on
$\left(\Sph, P\right)$ is the space of homeomorphisms
$\tau:\Sph\arr\CS$, where $\tau_1$ and $\tau_2$ are identified if
there exists a biholomorphic isomorphism $\Theta:\CS\arr\CS$
(i.e., an element of the M\"obius group) such that
$\Theta\circ\tau_1=\tau_2$ on $P$ and $\Theta\circ\tau_1$ is
isotopic to $\tau_2$ relative to $P$.

The \emph{moduli space} $\M_P$ of $\left(\Sph, P\right)$ is the
space of all injective maps $P\hookrightarrow\CS$ modulo
post-compositions with elements of the M\"obius group. The
Teichm\"uller space $\mathcal{T}_P$ is the universal cover of the
moduli space $\M_P$, where the covering map is the restriction map
of $\tau\in\mathcal{T}_P$ to $P$.

Let now $f:\Sph\arr \Sph$ be a branched covering with post-critical
set $P$. Then for every $\tau\in\mathcal{T}_P$ there exists a
unique element $\tau'\in\mathcal{T}_P$ such that we have a
commutative diagram
\begin{equation}\label{eq:pullback}\begin{array}{lcl}\Sph & \stackrel{f}{\arr} & \Sph\\
\mapdown{\tau'} & & \mapdown{\tau}\\
\CS & \stackrel{f_\tau}{\arr} & \CS \end{array}\end{equation} and
$f_\tau=(\tau')^{-1}\cdot f\cdot\tau:\CS\arr\CS$ is a rational
function. Let us write $\tau'=\sigma_f(\tau)$. The map
$\sigma_f:\mathcal{T}_P\arr\mathcal{T}_P$ is analytic and weakly
(i.e.\ non-uniformly) contracting (see~\cite{MR0302894}
or~\cite{DH:Thurston}).

Let us return to the case when $P$ is the post-critical set of the
``rabbit'' polynomial. The moduli space $\M_P$ is the set of maps
$\tau|_P:\{0, c, c^2+c, \infty\}\hookrightarrow\CS$ modulo
post-compositions with M\"obius transformations. We may assume,
applying an appropriate element of the M\"obius group, that $0$ is
mapped by $\tau$ to $0$, $c$ to $1$ and $\infty$ to $\infty$. Then
the points of moduli space are uniquely determined by the
value of $\tau|_P(c^2+c)=w$. We have $w\notin\{0, 1, \infty\}$,
since the map $\tau|_P$ is injective.

Therefore, the moduli space is isomorphic to the punctured sphere
$\CS\setminus\{0, 1, \infty\}$, where the point
$w\in\CS\setminus\{0, 1, \infty\}$ corresponds to the element
$\tau|_P$ such that $\tau|_P(\infty)=\infty, \tau|_P(0)=0,
\tau|_P(c)=1$ and $\tau|_P(c^2+c)=w$.

Let $f\in\mathcal{F}$ be arbitrary. Recall that $\mathcal{F}$ is the
set of degree-two branched coverings of $\CS$ with critical points $0,
\infty$ whose ramification graph coincides with that of the ``rabbit''
polynomial.

Let $\tau$ be arbitrary. Suppose that the projection of
$\tau'=\sigma_f(\tau)$ on moduli space is given by the point
$w_0\in\CS\setminus\{0, 1, \infty\}$ and the projection of $\tau$
is given by $w_1\in\CS\setminus\{0, 1, \infty\}$. Then the
rational function $f_\tau$ in the diagram~\eqref{eq:pullback} is a
degree-two map having critical points at $0$ and $\infty$ and
satisfying
\[f_\tau(\infty)=\infty,\quad f_\tau(0)=1,\quad f_\tau(1)=w_1, \quad f_\tau(w_0)=0,\]
since $f_\tau|_P=\left(\tau'|_P\right)^{-1}\cdot
f|_P\cdot\tau|_P$.

We conclude that $f_\tau$ is a quadratic polynomial. It is of the
form $az^2+1$, since $0$ is critical and $f_\tau(0)=1$. We get
therefore
\[a+1=w_1,\quad aw_0^2+1=0,\]
hence $a=-\frac{1}{w_0^2}$, so that
\[w_1=1-\frac{1}{w_0^2}.\]

We have thus obtained the following description of the action of the
pull-back map $\sigma_f$ on moduli space.

\begin{proposition}
\label{pr:sigmaf} The correspondence $\sigma_f(\tau)\mapsto\tau$
on Teichm\"uller space is projected on moduli space
$\M_P=\CS\setminus\{0, 1, \infty\}$ to the rational function
\[F:w\mapsto 1-\frac{1}{w^2}.\]
\end{proposition}

Suppose now that $h\in\mc$ is an arbitrary element of the mapping
class group of $\Ms=\CS\setminus P$. It defines an automorphism of
Teichm\"uller space by pre-composition: $\tau\mapsto
h\cdot\tau$. The mapping class group $\mc$ is the fundamental
group of the moduli space $\M_P$ and the action of $h$ on
$\mathcal{T}_P$ coincides with the corresponding deck
transformation. Note that if we identify elements of $\mathcal
T_P$ with paths in $\M_P$, then the action of $\mc$ by deck
transformations is given by pre-composition of paths; therefore
both actions of $\mc$ are left actions.

Let $f=f_R$ denote the ``rabbit'' polynomial. The corresponding point
of moduli space is given by the identical map $\{0, c, c^2+c,
\infty\}\arr\CS$. After normalization, we see that the
corresponding point of $\M_P=\C\setminus\{0, 1\}$ is the fixed
point $\frac{c^2+c}{c}=t\approx 0.8774 + 0.7449i$ of the rational
function $1-\frac{1}{w^2}$.

Let $\tau_0\in\mathcal{T}_P$ be the point of Teichm\"uller
space given by the identity map $\CS\arr\CS$. It is projected onto
the point $t$ of the moduli space $\M_P$. Every point
$\tau\in\mathcal{T}_P$ can be identified with the homotopy class
of a path $\ell_\tau$ in $\M_P$ starting at $t$, ending at the
projection of $\tau$ and equal to the image of a path in
$\mathcal{T}_P$ starting at $\tau_0$ and ending at $\tau$. The
homotopy class $\ell_\tau$ is uniquely defined and we have
$\ell_{h\cdot\tau}=\gamma_h\cdot\ell_\tau$, where
$\gamma_h\in\pi_1\left(\M_P, t\right)$ is the loop corresponding
to $h\in\mc$.

\begin{proposition}
\label{pr:sigmaell} For all $\tau\in\mathcal{T}_P$, $h\in\mc$ and
$f\in\mathcal{F}$ the following equalities hold:
\[\sigma_{h\cdot f}(\tau)=h\cdot\sigma_{f}(\tau),\quad
\sigma_{f\cdot h}(\tau)=\sigma_f(h\cdot\tau).\]

If $\tau\in\mathcal{F}$ corresponds to a path $\ell_\tau$ in
$\M_P$, then $\sigma_{f_R}(\tau)$ is represented by the path
$F^{-1}\left(\ell_\tau\right)\left[t\right]$, where
$F(w)=1-\frac{1}{w^2}$.
\end{proposition}

\begin{proof}
Consider the following commutative diagram:
\begin{equation}\label{eq:longdiagr}\begin{array}{lclclcl}\Sph & \stackrel{h}{\arr} & \Sph &
\stackrel{f}{\arr} & \Sph & \stackrel{h}{\arr} & \Sph\\
\mapdown{h\cdot\sigma_f(\tau)} & & \mapdown{\sigma_f(\tau)} & &
\mapdown{\tau} & &
\mapdown{h^{-1}\cdot\tau} \\
\CS & \stackrel{id}{\arr} & \CS & \stackrel{f_\tau}{\arr} & \CS &
\stackrel{id}{\arr} & \CS\end{array}\end{equation} It implies that
$\sigma_{h\cdot f}(\tau)=h\cdot\sigma_f(\tau)$ and $\sigma_{f\cdot
h}\left(h^{-1}\cdot\tau\right)=\sigma_f(\tau)$ for all $h\in\mc$
and $\tau\in\mathcal{T}_P$. The last equality implies that
$\sigma_{f\cdot h}(\tau)=\sigma_f(h\cdot\tau)$ for all
$\tau\in\mathcal{T}_P$.

The second statement of the proposition is a direct corollary of
Proposition~\ref{pr:sigmaf}.
\end{proof}

Consider now an arbitrary element $h\in\mc$ and the composition
$h\cdot f_R$. It is known (see~\cite{DH:Thurston}) that the orbit
of $\tau_0$ under iteration of $\sigma_{h\cdot f_R}$ will converge
to a point $\tau$ such that $\sigma_{h\cdot f_R}(\tau)=\tau$, and
the polynomial $f_\tau$ in the diagram~\eqref{eq:pullback} for
$f=h\cdot f_R$ is the polynomial which is Thurston equivalent to
$h\cdot f_R$.

Let $\gamma_h$ be the loop in $\M_P$ corresponding to $h\in\mc$.
By Proposition~\ref{pr:sigmaell}, we have $\sigma_{h\cdot
f_R}(\tau_0)=h\cdot\sigma_{f_R}(\tau_0)=h\cdot\tau_0$, hence the
path representing $\sigma_{h\cdot f_R}(\tau)$ is $\gamma_h$.

If $\ell_n$ is the path representing $\sigma_{h\cdot f_R}^{\circ
n}(\tau_0)$, then the path representing the point \[\sigma_{h\cdot
f_R}^{\circ (n+1)}(\tau_0)=h\cdot\sigma_{f_R}\left(\sigma_{h\cdot
f_R}^{\circ n}(\tau_0)\right)\] is $\gamma_h\cdot
F^{-1}\left(\ell_n\right)[t]$, by Proposition~\ref{pr:sigmaell}.
Consequently, the path representing the limit point $\tau$ is
\[\ell_\tau=\gamma_h^{(0)}\gamma_h^{(1)}\gamma_h^{(2)}\gamma_h^{(3)}\cdots,\]
where $\gamma_h^{(0)}=\gamma_h$ and $\gamma_h^{(n)}$ is the
preimage of $\gamma_h$ under $F^{\circ n}$ which starts at the end
of $\gamma_h^{(n-1)}$.

The endpoint of the path $\ell_\tau$ is one the three fixed points
of $1-\frac{1}{w^2}$. It is easy to see that the fixed point
$t\approx 0.8774 + 0.7449i$ corresponds to the ``rabbit'', the
point $\overline{t}\approx 0.8774 - 0.7449i$ corresponds to the
``corabbit'' and the point $\approx -0.7549$ corresponds to the
``airplane''.

\begin{figure}[ht]
\centering
  \includegraphics{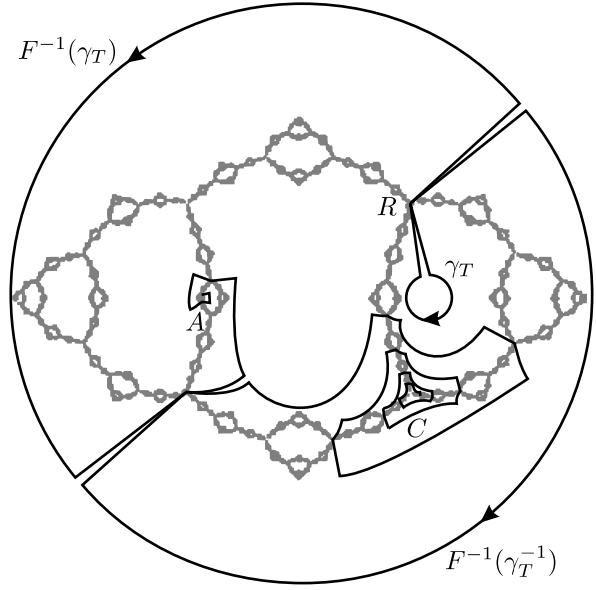}
  \caption{The path $\ell_\tau$ for $h=T$ and $T^{-1}$}\label{fig:twists}
\end{figure}

See for instance Figure~\ref{fig:twists}, where the paths
$\ell_\tau$ for $h=Y=T$ and $h=Y^{-1}=T^{-1}$ are indicated.

One can see that in the first case the path $\ell_\tau$ converges to
the fixed point corresponding to the ``airplane'' and in the second
case it converges to $\overline{t}$, which corresponds to the
``corabbit''. This can be shown after a detailed analysis of the
dynamics of the Fatou components of the rational function $F$.

Proposition~\ref{pr:psibar} can also be interpreted in these
terms: namely, we have $\psi=\phi_F$. If $h$ belongs to the domain
of the virtual endomorphism $\psi$, then the corresponding loop
$\gamma_h\in\pi_1\left(\M_P, t\right)$ belongs to the domain of
the virtual endomorphism $\phi_F$ associated with $F$. We have
then that the path converging to the fixed point of
$\sigma_{h\cdot f_R}$ is of the form
\[\gamma_h\cdot F^{-1}\left(\gamma_h\right)[t]
\cdot F^{-2}\left(\gamma_h\right)[t_1]\cdot
F^{-3}\left(\gamma_h\right)[t_2]\cdots,\] where $t_n$ is the
endpoint of $F^{-n}(\gamma_h)[t_{n-1}]$. This path is equal to
\[\gamma_h\gamma_g\cdot F^{-1}\left(\gamma_g\right)[t]
\cdot F^{-2}\left(\gamma_g\right)[t_1']\cdot
F^{-3}\left(\gamma_g\right)[t_2']\cdots,\] where
$g=\phi_F(h)=F^{-1}(h)[t]$ and $t_n'$ is the endpoint of
$F^{-n}\left(\gamma_g\right)[t_{n-1}']$. This proves that $h\cdot
f_R$ and $\psi(h)\cdot f_R$ are combinatorially equivalent.
Similar arguments work also in the case when $h$ does not belong
to the domain of $\psi$ (see~\cite[Theorem~6.6.3]{nek:book}).

\section{Preperiod $1$, period $2$}\label{se:3}

There are three families of quadratic topological polynomials with
three post-critical points: the first contains the ``rabbit'' and the
``airplane'', and its ramification graph is a cycle of length three.
The next family has ramification graph with preperiod $1$ and period
$2$; it contains the polynomial $f_i(z)=z^2+i$ and $f_{-i}=z^2-i$, as
well as obstructed topological polynomials. The last family has
ramification graph with preperiod $2$ and a fixed post-critical point;
it contains the polynomials $\approx z^2-1.5434$ and $\approx
z^2-0.2282\pm 1.1151 i$, and is dealt with in Section~\ref{s:quater}.

\subsection{The iterated monodromy group}
Let us consider first the polynomial $f_i(z)=z^2+i$. The dynamics of
$f_i$ on its post-critical set is
\[i\mapsto -1+i\mapsto -i\mapsto -1+i.\]

\begin{figure}[ht]
\centering
  \includegraphics{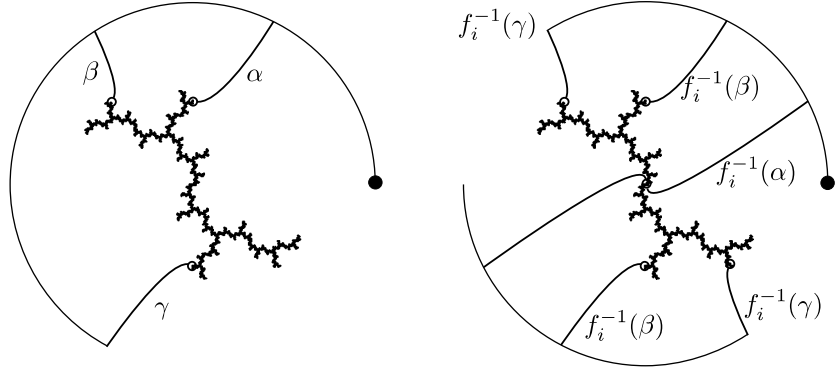}
  \caption{Computation of $\img{z^2+i}$}\label{fig:i}
\end{figure}

Let us compute the iterated monodromy group of $f_i$. We again
take $+\infty$ as the basepoint. Let $\alpha$ be the loop going
around $i$ in the positive direction and connected to the
basepoint by the external ray $R_{1/6}$ and the arc $[0, 1/6]$ of
the circle at infinity. The loops $\beta$ and $\gamma$ go around
the points $-1+i$ and $-i$ and are connected to the circle at
infinity by the rays $R_{1/3}$ and $R_{2/3}$ and to the point
$+\infty$ by the arcs $[0, 1/3]$ and $[0, 2/3]$, respectively. The
connecting paths $\ell_0$ and $\ell_1$ are, as usual, the trivial
path and the upper semicircle. See the loops $\alpha, \beta,
\gamma$ and their preimages in Figure~\ref{fig:i}. Computation of
the wreath recursion gives
\[\Phi_{f_i}(\alpha)=\pair<\alpha^{-1}\beta^{-1}, \beta\alpha>\sigma,
\quad\Phi_{f_i}(\beta)=\pair<\alpha,
\gamma>,\quad\Phi_{f_i}(\gamma)=\pair<\beta, 1>.\]

We see that the corresponding elements $\alpha, \beta, \gamma$ of
$\img{f_i}$ are of order two.

Let $\F$ be the set of homotopy classes of branched coverings which
have the same ramification graph as $f_i$. Namely, $\F$ is the set of
homotopy classes of degree-two topological polynomials $f$ such that
the finite post-critical set of $f$ has three different points $c_1,
c_2, c_3$ which are mapped in the following way:
\[c_1\mapsto c_2\mapsto c_3\mapsto c_2\]
(where we assume that the same points $c_1, c_2, c_3$ are chosen
for all elements of $\F$).

Among quadratic polynomials $f(z)=z^2+c$ only $z^2+i$ and $z^2-i$
have this ramification graph, since the set of the roots of
the equation $f^3(c)=f(c)$ is $\{0, -1, i, -i, -2\}$, but $z^2,
z^2-1$ and $z^2-2$ have different post-critical dynamics.

\begin{figure}[ht]
\centering
\includegraphics{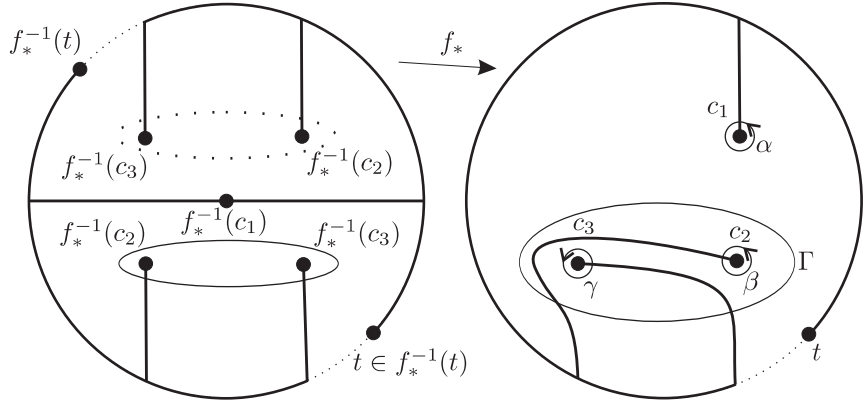}
  \caption{An obstructed topological polynomial}\label{fig:obstr}
\end{figure}

Note, however, that there exist obstructed topological polynomials
in $\F$. A way to construct one is shown in Figure~\ref{fig:obstr}.
It describes the post-critical points $c_1, c_2, c_3$ and curves
connecting them to infinity on the right-hand side of the figure.
The left-hand side shows the preimages of the points and curves.
The map folds the horizontal line in two and maps the critical
point $f_*^{-1}(c_1)$ to $c_1$. It is a homeomorphism from each of
the upper and lower half-planes to the complement of the line
connecting $c_1$ to infinity.

Consider the simple closed curve $\Gamma$ around the points $c_2$
and $c_3$. It has two $f_*$-preimages. One is peripheral, and the
other is homotopic to $\Gamma$. The map $f_*$ is of degree $1$ on
the non-peripheral preimage of $\Gamma$, so the curve $\Gamma$ is
an obstruction.

Let us compute the iterated monodromy group of the map $f_*$.  Choose
the basepoint $t$ on the circle at infinity, as shown on
Figure~\ref{fig:obstr}, and let the connecting paths from $t$ to its
preimages be the trivial path $\ell_0$ at $t$ and the semicircle
$\ell_1$ going in the positive direction from $t$ to its other
preimage. Let $\alpha$, $\beta$, $\gamma$ be generators of
$\pi_1(\CS\setminus P_{f_*},t)$ following the circle at infinity and
the arcs towards the points $c_1, c_2, c_3$ respectively, encircling
the point in the positive direction, and returning back to $t$, as
shown on the right-hand side of Figure~\ref{fig:obstr}.

It is easy to see from Figure~\ref{fig:obstr} that the wreath
recursion is
\begin{equation}\label{eq:fstar}\Phi_{f_*}(\alpha)=\pair<\alpha^{-1},
\alpha>\sigma, \quad\Phi_{f_*}(\beta)=\pair<\alpha, \gamma>,\quad
\Phi_{f_*}(\gamma)=\pair<1, \gamma\beta\gamma^{-1}>.\end{equation}

Direct computation shows that the corresponding elements $\alpha,
\beta, \gamma$ of the iterated monodromy group are of order two
and that $\beta$ and $\gamma$ commute. Conjugating $\img{f_*}$ by
$\Delta=\pair<\alpha\Delta, \Delta>$, we see that the generators
of $\img{f_*}$ can be defined by the recursion
\[\overline\alpha=\sigma,\quad\overline\beta=\pair<\overline\alpha,\overline\gamma>,\quad
\overline\gamma=\pair<1, \overline\beta>.\]

This is one of the Grigorchuk groups $G_\omega$
from~\cite{grigorchuk:growth_en}, namely that given by
$\omega=(01)^\infty$.  One of the ways to see that it is not the
iterated monodromy group of a rational function (and thus that the map
$f_*$ is not combinatorially equivalent to a polynomial) is to show
that its limit orbispace (see~\cite[Section~4.6]{nek:book}) has an
isotropy group isomorphic to the Klein group $C_2\times C_2$, namely
$\langle\overline\beta,\overline\gamma\rangle$.

Note also that the composition of $f_*$ with a power of the Dehn twist
about the curve $\Gamma$ is also obstructed (with the same
obstruction $\Gamma$).

\subsection{Moduli space}
Let us compute the virtual endomorphism $\psi$ on the mapping
class group $\mc$ of the punctured plane $\C\setminus
P=\C\setminus\{i, i-1, -i\}$, via the moduli space approach.

Let $\F$ be as before the set of homotopy classes of branched
coverings having the same ramification graph as $z^2+i$. The moduli
space of $\CS\setminus P$ is again a plane with two punctures. We
normalize the mappings $P\hookrightarrow\CS$ so that $\infty$ is mapped to $\infty$,
$0$ is the critical value of $f\in \F$ and $f(0)=1$. Let us compute the action of
the inverse of $\sigma_f$ on moduli space in this case. We have
\[f_\tau(0)=1,\quad f_\tau(1)=w_1,\quad f_\tau(w_0)=1.\]
We also have that $0$ is the critical value of the quadratic
polynomial $f_\tau$. Hence $f_\tau$ is of the form $(az-1)^2$,
where $1/a$ is its critical point. The last equality implies that
$aw_0-1=1$ (as $aw_0-1=-1$ would give $w_0=0$), hence $a=2/w_0$, so
\[w_1=f_\tau(1)=\left(\frac{2-w_0}{w_0}\right)^2.\]

Therefore, the correspondence $\sigma_f(\tau)\mapsto\tau$ is
projected to the rational function
$F(z)=\left(\frac{2-w}{w}\right)^2$ on moduli space (compare
with Proposition~\ref{pr:sigmaf}).

The fixed points of the rational function
$\left(\frac{2-w}{w}\right)^2$ are $w=1$ and $w=\pm 2i$. The point $1$
belongs to the post-critical set and is the puncture of the moduli
space. Thus, it does not correspond to any quadratic polynomial (and,
as we will see, corresponds to obstructed topological polynomials).
The point $2i$ corresponds to the polynomial
\[\frac{4}{(2i)^2}\left(z-\frac{2i}{2}\right)^2=-\left(z-i\right)^2,\]
which is conjugate to $z^2+i$. The point $-2i$ corresponds to
$-\left(z+i\right)^2$, which is conjugate to $z^2-i$.

The critical points of the rational function
$\left(\frac{2-w}{w}\right)^2$ are $w=2, w=0$. We have the
following dynamics on the post-critical orbit:
\[2\mapsto 0\mapsto \infty\mapsto 1\mapsto 1.\]

Therefore, the post-critical set of
$F(w)=\left(\frac{2-w}{w}\right)^2$ is $\{0, 1, \infty\}$. The
Julia set of $F$ is the whole sphere, since there are no
superattracting cycles. Actually, the Thurston orbifold
$\mathcal{O}_F$ of $F$ is $(2, 4, 4)$, i.e., the orbifold of the
action on $\C$ of the group of affine transformations $\{z\mapsto
i^kz+a+bi\;:\;k, a, b\in\Z\}$. Moreover, $F$ is induced on the
orbifold $\mathcal{O}_F$ by the expanding map $z\mapsto (1+i)z$ of
$\C$.

To show this explicitly, consider the Weierstrass function
\[\wp(z)=\frac{1}{z^2}+\sum_{\omega\in\Z[i]\setminus\{0\}}\left[\frac{1}{(z+\omega)^2}-\frac{1}{\omega^2}\right]\]
associated to the lattice $\Z[i]$.

It follows from its definition that $\wp$ is an even function and
that, by the choice of the lattice $\Z[i]$, we have $\wp(iz)=-\wp(z)$.
Consequently, $i\wp'(iz)=-\wp'(z)$ and thus
\[\wp'(iz)=i\wp'(z).\]

It is known that
\[(\wp'(z))^2=4\wp^3(z)-g_2\wp(z)-g_3,\] with $g_2=60s_4$ and
$g_3=140s_6$, where $s_m=\sum_{\omega\in\Z[i], \omega\ne 0}
\omega^{-m}$ are the Eisenstein series.  It is clear that $s_6=0$
because $\Z[i]$ has a $4$-fold symmetry, so $g_3=0$.

Another classical fact is the addition formula
\[\wp(z_1+z_2)=-\wp(z_1)-\wp(z_2)+
\frac{1}{4}\left(\frac{\wp'(z_1)-\wp'(z_2)}{\wp(z_1)-\wp(z_2)}\right)^2.\]

Let us then compute $\wp((i+1)z)$:
\begin{multline*}
  \wp(iz+z)=-\wp(iz)-\wp(z)+
\frac{1}{4}\left(\frac{\wp'(iz)-\wp'(z)}{\wp(iz)-\wp(z)}\right)^2=
\frac{1}{4}\left(\frac{(i-1)\wp'(z)}{-2\wp(z)}\right)^2\\
=-\frac{i}{8}\cdot\frac{4\wp^2(z)-g_2}{\wp(z)}.
\end{multline*}

The function $\wp:\C\arr\CS$ realizes the universal covering of
the orbifold $(2, 2, 2, 2)$. The group of deck transformations of
this covering is the group of all affine transformations of the
form $z\mapsto \pm z+c$, with $c\in\Z[i]$. In other words this
group is the group of holomorphic automorphisms $\alpha:\C\arr\C$
such that $\wp(\alpha(z))=\wp(z)$.

It follows that the function\footnote{Here
  $\wp^2(z)=\left(\wp(z)\right)^2$.} $\wp^2:\C\arr\CS$ realizes the
universal covering of the orbifold $(2, 4, 4)$ for which the group
$G$ of deck transformations is the group of affine transformations
of the form $z\mapsto i^kz+c$, with $c\in\Z[i]$.

We have
\[\wp^2\left((i+1)z\right)=-\frac{1}{64}\left(\frac{4\wp^2(z)-g_2}{\wp(z)}\right)^2=
-\frac{\left(\wp^2(z)-\frac{g_2}{4}\right)^2}{4\wp^2(z)}.\]

It follows that the map $z\mapsto (i+1)z$ on the universal cover
$\C$ induces on the orbifold $G\backslash\C$ the rational map
$t\mapsto-\frac{\left(t-\frac{g_2}{4}\right)^2}{4t}$. Performing
the change of variables $t=\frac{g_2/4}{1-w}$ we see that this
rational map is conjugate to
$w\mapsto\left(\frac{2-w}{w}\right)^2$.

We conclude that the iterated monodromy group of the rational
function $\left(\frac{2-w}{w}\right)^2$ is isomorphic to the group
of affine transformations $z\mapsto i^kz+c$, $c\in\Z[i]$
(see~\cite[Subsection~6.3.2.2]{nek:book}
and~\cite[Section~5]{bgn}).

\begin{figure}[ht]
\centering
  \includegraphics{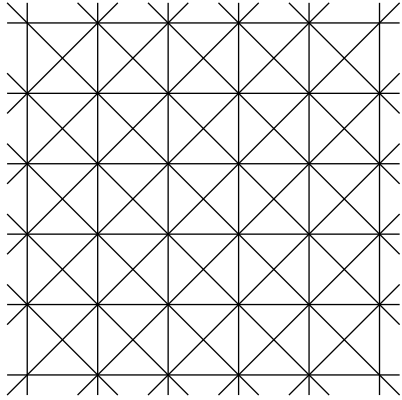}
  \caption{The fundamental group of the orbifold $(2, 4, 4)$}\label{fig:lattice1}
\end{figure}

This group is the group of the orientation-preserving automorphisms of
the tiling of the plane by triangles shown on
Figure~\ref{fig:lattice1}. The triangles are fundamental domains of
the group action. We assume that the vertices of the grid coincide
with the Gaussian integers $\Z[i]$.

Let us cut the Riemann sphere $\CS$ along the ray $[0,+\infty]$
consisting of the non-negative reals and infinity. This ray will
then contain the post-critical set of
$F(z)=\left(\frac{2-z}{z}\right)^2$. It is also easy to see that
the preimage of this cut in the universal cover $\C$ of the
orbifold $\mathcal{O}_F$ is precisely the union of the lines of
the tiling in Figure~\ref{fig:lattice1}. In particular, the
preimage of $\CS\setminus[0,\infty]$ is the disjoint union of the
open triangles of the tiling.

\subsection{The iterated monodromy group on moduli space}\label{ss:affine}
\begin{figure}[ht]
\centering
  \includegraphics{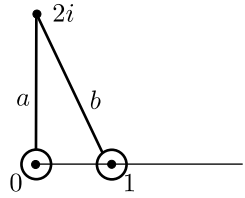}
  \caption{The generators of $\img{F}$}\label{fig:ab}
\end{figure}

\begin{figure}[ht]
\centering\includegraphics{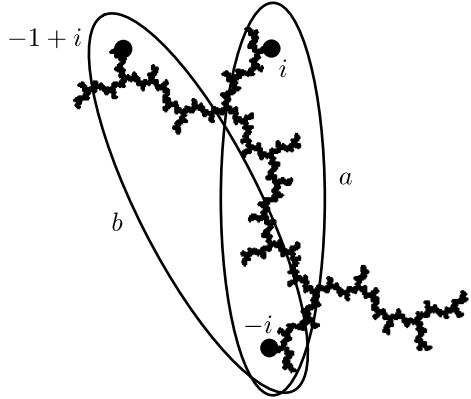} \caption{The generators $a$
and $b$ as Dehn twists}\label{fig:abt}
\end{figure}

Let us compute the recursion defining $\img{F}$. We take $t=2i$ as
basepoint. The iterated monodromy group is generated by the
loops $a$ and $b$ going in the positive direction around $0$ and
$1$ respectively, as shown in Figure~\ref{fig:ab}. The loops $a$
and $b$ correspond to the right Dehn twists about the curves going
around $i$, $-i$ and $-1+i, -i$, respectively (see
Figure~\ref{fig:abt}).

\begin{figure}[ht]
\centering
  \includegraphics{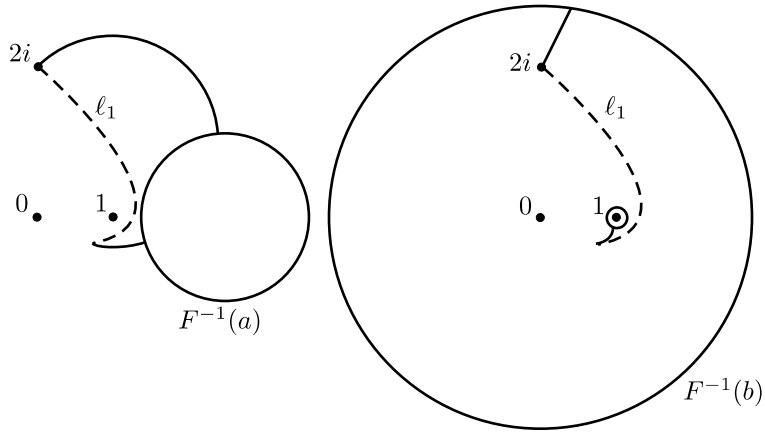}
  \caption{Computation of $\img{F}$}\label{fig:imod}
\end{figure}

The preimages of $2i$ under $F$ are $2i$ and $\frac{4-2i}{5}$. Let
$\ell_0$ be the trivial path at $2i$ and let $\ell_1$ be the path
connecting $2i$ to $\frac{4-2i}{5}$ passing between $1$ and $2$ as
shown in Figure~\ref{fig:imod}.

The preimages of the loops $a$ and $b$ are shown on
Figure~\ref{fig:imod}. We see that the wreath recursion is
\begin{equation}\label{eq:abrec}\Phi(a)=\sigma,\qquad\Phi(b)=\pair<b^{-1}a^{-1},
b>.\end{equation}

The virtual endomorphism associated with the first coordinate of
the wreath recursion $\Phi$ on the mapping class group is given by
\begin{equation}\label{eq:psiab}\psi(a^2)=1,\quad\psi(b)=b^{-1}a^{-1},\quad\psi(b^a)=b,\end{equation}
since
\[\Phi(a^2)=1,\quad\Phi(b)=\pair<b^{-1}a^{-1}, b>,\quad
\Phi(b^a)=\pair<b, b^{-1}a^{-1}>.\]

This virtual endomorphism has the property that $f_i\cdot
g=\psi(g)\cdot f_i$ for every $g$ in the index-two subgroup $\dom\psi=\langle
a^2, b, b^a\rangle$ of the mapping class group $\mc$. The proof is the same as in
the case of the ``rabbit'' polynomial.

Let us see now what happens on the universal cover of the orbifold
$\mathcal{O}_F$. For every path $\gamma$ in the moduli space
$\M_P$ and for every preimage $\zeta\in\C$ of the startpoint of
$\gamma$ in the universal cover of the orbifold $\mathcal{O}_F$,
there exists a unique preimage of $\gamma$ in $\C$ which starts at
$\zeta$.

Let $t\in\C$ be the preimage of $2i$ under the covering map
$\C\arr\mathcal{O}_F$ which belongs to the triangle $\Delta$ with
vertices $0, 1, (1+i)/2$. Note that $0$ and $1\in\C$ are preimages of
$1\in\mathcal{O}_F$, that $(1+i)/2$ is a preimage of
$\infty\in\mathcal{O}_F$, and that $1/2$ is a preimage of
$0\in\mathcal{O}_F$.

\begin{figure}[ht]
\centering
  \includegraphics{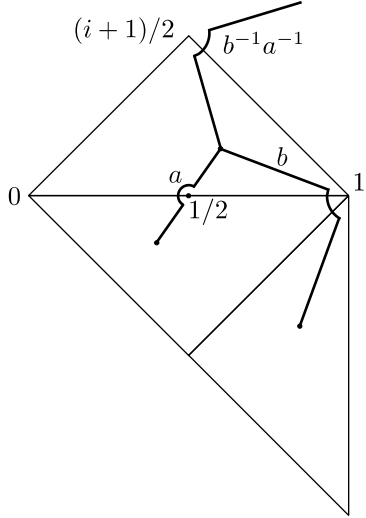}
  \caption{Loops $a$ and $b$ in the universal cover $\C$}\label{fig:ab2}
\end{figure}

Then the preimages of the loops $a$ and $b$ starting at $t$ go in
the positive direction around the points $1/2$ and $1$, if we assume
that the universal covering map $\C\arr\mathcal{O}_F$ is
orientation-preserving (see Figure~\ref{fig:ab2}). We use here
the fact that the segment $[1/2, 1]$ is a preimage of the interval $[0,
1]\subset\mathcal{O}_F$. It follows that the elements $a$ and $b$
of the fundamental group of $\M_P$ act on the plane $\C$ via the
affine maps
\[A(z)=-z+1,\quad B(z)=iz+1-i,\]
respectively. Recall that this has to be a left action.

The element $b^{-1}a^{-1}$ of the fundamental group $\pi_1(\M_P,
2i)$ is a small loop around $\infty$ connected to the basepoint by
a curve disjoint from the positive ray. It follows that the
preimage of the loop $b^{-1}a^{-1}$ in the triangle $\Delta$  is
the top curve shown in Figure~\ref{fig:ab2}.

We have \[B^{-1}A^{-1}(z)=iz+1,\quad B^A(z)=iz.\] The preimage of
the loop $b^a$ in the triangle $\Delta$ goes around the vertex
$0$.

Let $\mathcal{T}_P$ be the Teichm\"uller space modeled on $(\CS,
P)$ for $P=\{i, -1+i, -i, \infty\}$. We have the following
commutative diagram
\[\begin{array}{lcl}\mathcal{T}_P & \arr & \C\\
\downarrow & & \downarrow\\
\M_P & \hookrightarrow & \mathcal{O}_F \end{array}\] where
$\mathcal{T}_P\arr\M_P$ is the natural projection,
$\M_P\hookrightarrow\mathcal{O}_F=\CS$ is the identical
embedding and $\C\arr\mathcal{O}_F$ is the universal covering
map. The map $\mathcal{T}_P\arr\C$ exists, since $\mathcal{T}_P$
is the universal cover of $\M_P$. It is defined uniquely up to
an element of the group $G=\{i^kz+c\;:\;c\in\Z[i]\}$ acting on
$\C$.

Let $\tau_0\in\mathcal{T}_P$ be the original complex structure on
$\CS$. It is projected to $2i\in\M_P$ and we may assume that its
image $\zeta$ in $\C$ belongs to the triangle $\Delta$. This fixes
uniquely the map $\mathcal{T}_P\arr\C$.

We know that the correspondence $\tau\mapsto\sigma_{f_i}(\tau)$ on
$\mathcal{T}_P$ is semi-conjugated via the map $\mathcal{T}_P\arr\C$
to an affine map $\Sigma:z\mapsto i^k\frac{z}{1+i}+\xi$ for some
$k\in\Z$ and $\xi\in\C$ (since $F$ is induced by multiplication by
$(1+i)$ on $\C$). The point $\zeta$ is fixed by $\Sigma$, since
$\tau_0$ is fixed under $\sigma_{f_i}$.

\begin{figure}[ht]
\centering
  \includegraphics{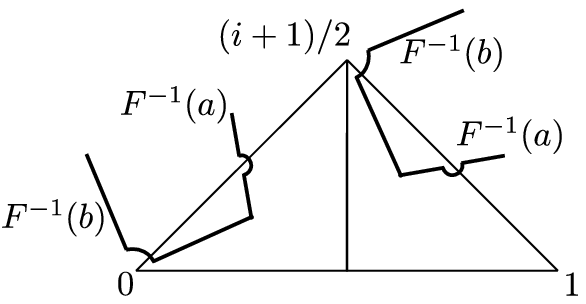}
  \caption{}
  \label{fig:ab3}
\end{figure}

The preimage of the positive ray under the action of $F$ is the whole
real line. The preimage of the negative ray in the universal cover
$\C$ divides the triangle $\Delta$ into two similar triangles, see
Figure~\ref{fig:ab3}. This figure also displays the preimages of the
curves $F^{-1}(a)$ and $F^{-1}(b)$.

We know from the wreath recursions~\eqref{eq:abrec} that the
$F$-preimage of $b$ which starts at $2i$ coincides with
$b^{-1}a^{-1}$, which is a loop around infinity. This shows that
the point $\zeta$ belongs to the right-hand side triangle on
Figure~\ref{fig:ab3}.

It follows that \[\Sigma(z)=\frac{i-1}{2}z+1\text{ and, therefore,
}\zeta=\frac{1+i}{2+i}=\frac{3+i}{5}.\]

\subsection{The answer}
We are now in a position to determine the combinatorial
equivalence class of the Thurston maps considered in the
non-obstructed case. The obstructed maps will be considered later.

Let $\pi$ be the homomorphism from the mapping class group $\mc$
onto the affine group $G=\{z\mapsto i^kz+c\;:\;k\in\Z,
c\in\Z[i]\}$ defined on the generators by
\[\pi(a)=A:z\mapsto -z+1,\quad \pi(b)=B:z\mapsto iz+1-i.\]
We have shown above that $\pi$ is the canonical epimorphism of $\mc$
onto the iterated monodromy group
$\img{\left(\frac{2-z}{z}\right)^2}$. Note that multiplication in the
affine group $G$ corresponds to left action, i.e.,
$\pi(g_1g_2)=\pi(g_1)\circ\pi(g_2)$.

\begin{figure}[ht]
\centering
  \includegraphics{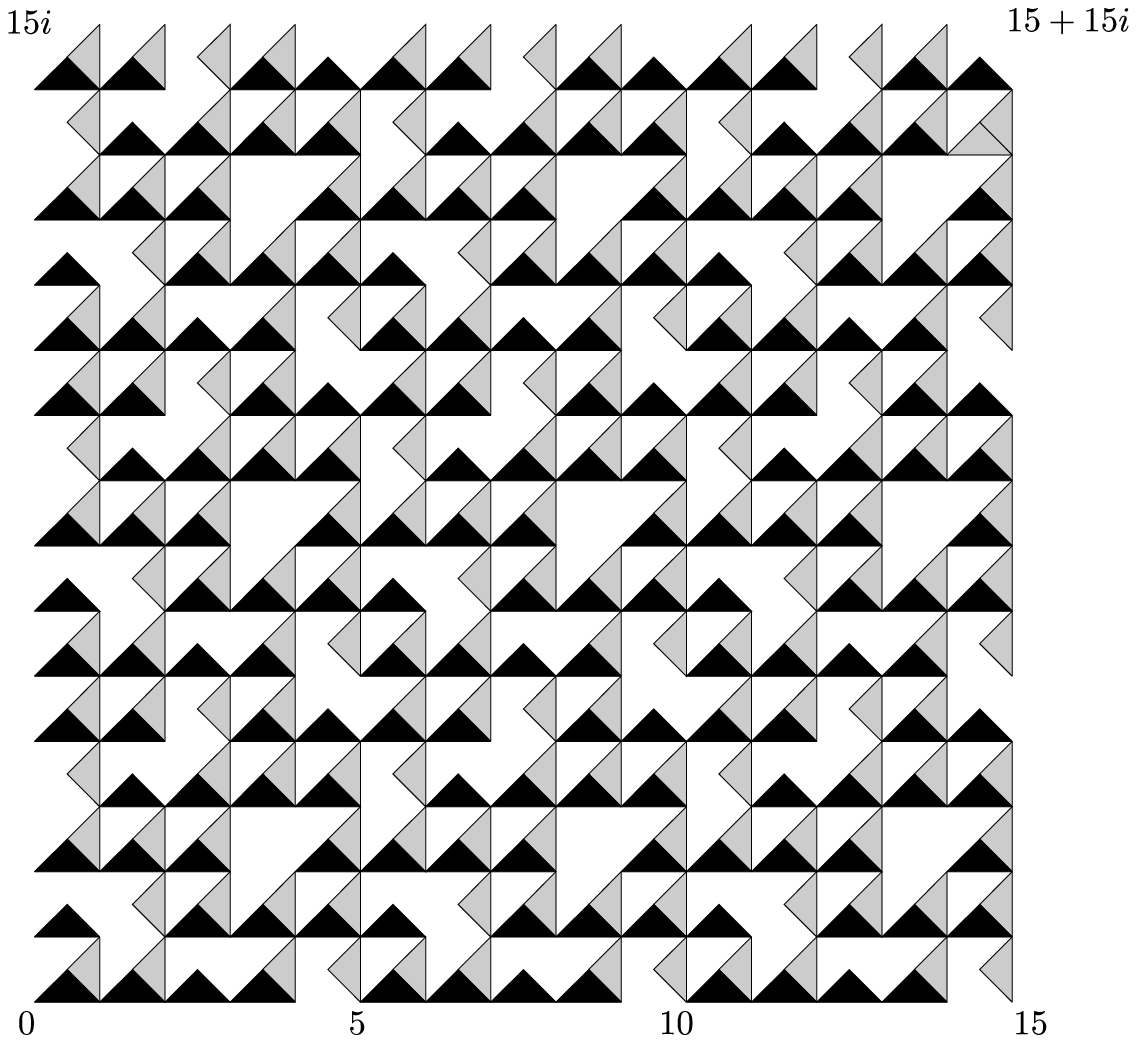}
  \caption{}\label{fig:lattice3}
\end{figure}

\begin{theorem}\label{th:mod5}
  Let $h\in\mc$ be an element of the mapping class group.
  Write $\pi(h)$ as the affine transformation $\pi(h)(z)=i^kz+c$. Then:
  \begin{itemize}
  \item if $k=0$ and $\frac{c-i-1}{2+i}\notin\Z[i]$,
    then $f_i\cdot h$ is equivalent to $f_i$;
  \item if $k=1$ and $\frac{c-i-1}{1+2i}\notin\Z[i]$,
    then $f\cdot h$ is equivalent to $f_{-i}$;
  \item in all other cases, $f$ is obstructed.
  \end{itemize}
  Let $Q$ be the group $(\Z[i]/5\Z[i])\rtimes\Z/4$ of order 100, where
  the action of $\Z/4$ is by multiplication by $i$.  Then the answer
  (in $\{f_i,f_{-i},\text{obstructed}\}$) depends only on the image of
  $h$ in $Q$ under the homomorphism $\mc\arr Q$ mapping $a$ to $(-1,2)$
  and $b$ to $(1-i,1)$.
\end{theorem}
Therefore, unlike in the ``rabbit/airplane'' case, the answer is
periodic.  See Figure~\ref{fig:lattice3} for dependence of the answer
on the element $\pi(h)$ of the group $G$. The lower left black
triangle is the origin (it corresponds to the trivial element of $G$).
For every other triangle $\Delta$ of the picture there is a unique
element $h$ mapping the original triangle to $\Delta$. If $\Delta$ is
black, then $f_i\cdot h$ is equivalent to $z^2+i$, if $\Delta$ is
grey, then $f_i\cdot h$ is equivalent to $z^2-i$. White triangles
correspond to obstructed maps. The picture is periodic with period $5$
in both directions.

\begin{proof}[Proof of Theorem~\ref{th:mod5}]
  By Proposition~\ref{pr:sigmaell}, the map $\sigma_{f_i\cdot h}$ is
  projected to the affine transformation $\Sigma\circ\pi(h)$ of $\C$.  The
  fixed point of $\sigma_{f_i\cdot h}$, if it exists, is mapped to the
  fixed point of $\Sigma\circ H(z)=\frac{i-1}{2}(i^kz+c)+1$ in $\C$,
  which is
  \[\xi=\frac{\frac{i-1}{2}c+1}{1-\frac{i-1}{2}i^k}=-\frac{c-i-1}{i+1+i^k}.\]

  Let us consider the possible values for $k$ in turn.
  \begin{description}
  \item[if $\boldsymbol{k=0}$] then $\xi=-\frac{c-i-1}{i+2}$, which is
    in the same $G$-orbit as $\frac{c+1}{i+2}$. The possible residues
    modulo $i+2$ are $0, \pm i, \pm 1$, which, modulo multiplication
    by powers of $i$ are either $0$ or $1$. Consequently, if
    $\frac{c+1}{i+2}\in\Z[i]$, then $\xi$ belongs to the $G$-orbit of
    $0$, otherwise it belongs to the $G$-orbit of $\frac{1}{i+2}$,
    which coincides with the $G$-orbit of $\zeta=\frac{3+i}{5}$.

    If $\xi$ belongs to the orbit of $0$, then the map $f_i\cdot h$ is
    obstructed. If it belongs to the orbit of $\zeta$, then $f_i\cdot h$
    is equivalent to $f_i(z)=z^2+i$.

  \item[if $\boldsymbol{k=1}$] then $\xi=-\frac{c-i-1}{2i+1}$, which
    is in the same $G$-orbit as $\frac{c+i}{2i+1}$. Here again, if
    $\frac{c+i}{2i+1}\in\Z[i]$, then $f_i\cdot h$ is obstructed,
    otherwise it is equivalent to $z^2-i$.

  \item[if $\boldsymbol{k=2}$] then $\xi=-\frac{c-i-1}{i}$, which is
    integral, and we get an obstructed map.

  \item[if $\boldsymbol{k=3}$] then $\xi=-c+i+1$, and we also get an
    obstructed map.
  \end{description}
\end{proof}

\subsection{The obstructed cases} No description of combinatorial
classes of obstructed Thurston maps follows directly from Thurston's
Theorem~\ref{th:thurston}; therefore we have to treat the obstructed
cases separately.

Let $\mathcal{F}$ be the set of quadratic topological polynomials
$f:\CS\arr\CS$ with finite critical value $i$ and such that
$f(i)=-1+i$, $f(-1+i)=-i$ and $f(-i)=-1+i$. Let $\mathfrak{F}$ be
the set of homotopy classes of elements of $\mathcal{F}$ relative
to $P=\{\infty, i, -1+i, -i\}$. Let $\mc$ be the mapping class
group of $\CS\setminus P$. Recall that $\F$ is a $\mc$-bimodule by
pre- and post-composition.

\begin{proposition}\label{prop:irred}
For every $f\in\mathfrak{F}$ there exist elements $h$ and
$g\in\mc$ such that $h\cdot f\cdot g$ is equal to the polynomial
$f_i(z)=z^2+i$.
\end{proposition}

In the terminology of~\cite{nek:book} this means that the
$\mc$-bimodule $\mathfrak{F}$ is \emph{irreducible}.

\begin{proof}
Let $\mathcal{T}_P$ and $\M_P$ denote the Teichm\"uller and the moduli
space of $\CS\setminus P$, respectively. We identify, as above,
$\M_P$ with $\CS\setminus\{0, 1, \infty\}$ so that the
correspondence $\sigma_f(\tau)\mapsto\tau$ on $\mathcal{T}_P$ is
projected to the rational function
$F(w)=\left(\frac{2-w}{w}\right)^2$ on $\M_P$.

Let $\zeta=2i$ be the fixed point of $F$ corresponding to the
polynomial $z^2+i$. Let $\tau_i$ be an arbitrary preimage of
$\zeta$ in $\mathcal{T}_P$. Then the image of $\sigma_f(\tau_i)$
in $\M_P$ is equal to one of the $F$-preimages of $\zeta$, i.e.,
either to $2i$ or to $\frac{4-2i}{5}$.

The path $\sigma_{f}(a)$ is a lift to $\mathcal{T}_P$ of an
$F$-preimage of the loop $a$ in $\M_P$. We know, by~\eqref{eq:abrec},
that both $F$-preimages of $a$ are paths starting at one of the points
$2i$, $\frac{4-2i}{5}$ and ending at the other.

Consequently, if the image of $\sigma_f(\tau_i)$ in $\M_P$ is
equal to $\frac{4-2i}{5}$, then the image of $\sigma_{f\cdot
a}(\tau_i)=\sigma_{f}(a\cdot\tau_i)$ is equal to $2i$ (see
Proposition~\ref{pr:sigmaell}).

Hence, either $\sigma_f(\tau_i)$ or $\sigma_{f\cdot a}(\tau_i)$
belongs to the $\pi_1(\M_P)$-orbit of $\tau_i$, i.e., there exists
a homeomorphism $h\in\mc$ such that $\sigma_{h\cdot
f}(\tau_i)=h\cdot\sigma_f(\tau_i)=\tau_i$, or $\sigma_{h\cdot
f\cdot a}(\tau_i)=h\cdot\sigma_{f\cdot a}(\tau_i)=\tau_i$. But
this implies that either $h\cdot f$ or $h\cdot f\cdot a$ is
combinatorially equivalent to $f_i$.
\end{proof}

\begin{corollary}
\label{cor:lefttransitive} For every branched covering
$f\in\mathfrak{F}$ there exists $g\in\mc$ such that $f=f_i\cdot
g$.
\end{corollary}

\begin{proof}
  This follows from Proposition~\ref{prop:irred} and the fact that the
  associated virtual endomorphism $\psi$ is onto,
  see~\eqref{eq:psiab}.
\end{proof}


It will be more convenient for us to proceed with the branched
covering $f_i\cdot a$, which is obstructed by Theorem~\ref{th:mod5}.
We now show that it is equivalent to the branched covering $f_*$ shown
on Figure~\ref{fig:obstr}.

For that purpose, note that the twist $a$ acts on the generators
$\alpha, \beta, \gamma$ of the fundamental group of the punctured
plane in the same way as the element $h$ in~\eqref{eq:h123}.
Applying that twist to the wreath recursion $\Phi_{f_i}$, we see
that the wreath recursion $\Phi_{f_ia}$ is given
by~\eqref{eq:fstar}, so that the maps $f_ia$ and $f_*$ are
equivalent, by Proposition~\ref{pr:homotopy}.

Let $\phi$ denote the virtual endomorphism of $\mc$ such that
$f_*\cdot g=\phi(g)\cdot f_*$. It follows from the wreath
recursion $\Phi$ defined by~\eqref{eq:abrec} that $\phi$ is the
virtual endomorphism associated with the second coordinate of the
wreath recursion $\Phi$ and thus has domain $H=\langle
a^2,b,b^a\rangle$ and is given by
\begin{equation}\label{eq:phi}\phi(a^2)=1,\quad\phi(b)=b,\quad\phi(b^a)=b^{-1}a^{-1}.\end{equation}

Let us introduce the following subgroups $\mathcal{E}_n\lhd\mc$
(see~\cite[Subsection~2.13.1]{nek:book}). We set
$\mathcal{E}_0=\{1\}$ and
\[\mathcal{E}_n=\{g\in H<\mc\;:\;\text{$\phi(g)$ and $\phi(g^a)$
belong to $\mathcal{E}_{n-1}$}\}.\] In other words,
$\mathcal{E}_n$ is the kernel of the wreath recursion
$\Phi^n:\mc\arr\mc\wr\symm[X^n]$ describing the action and the
restrictions of $\mc$ on words of length $n$.

It follows that $\mathcal{E}_n$ are normal subgroups and that
$\mathcal{E}_n\ge\mathcal{E}_{n-1}$. Let us denote
$\Einfty=\bigcup_{n\ge 0}\mathcal{E}_n$.

\begin{lemma}
\label{l:en} If $g\in\Einfty$, then for all
$f\in\mathfrak{F}$ the branched coverings $f\cdot g$ and $f$ are
combinatorially equivalent.
\end{lemma}

\begin{proof}
We have $f=f_*\cdot h$ for some $h\in\mc$, by
Corollary~\ref{cor:lefttransitive}. Let us prove the statement by
induction on $n$ such that $g\in\mathcal{E}_n$. The statement is
trivial for $n=0$. We have
\begin{multline*}f_*\cdot
hg=f_*\cdot
hgh^{-1}\cdot h\\
= \phi(hgh^{-1})\cdot f_*\cdot h=\phi(hgh^{-1})\cdot\left(f_*\cdot
h\phi(hgh^{-1})\right)\cdot\phi(hgh^{-1})^{-1}.
\end{multline*}
We have $hgh^{-1}\in\mathcal{E}_n$, hence
$\phi(hgh^{-1})\in\mathcal{E}_{n-1}$, and by the inductive
hypothesis, $f_*\cdot h\phi(hgh^{-1})$ is combinatorially
equivalent to $f_*\cdot h$, which shows that $f_*\cdot hg$ is
combinatorially equivalent to $f_*\cdot h$.
\end{proof}

\begin{corollary}
\label{cor:einfty} The combinatorial equivalence class of a
branched covering $f_*\cdot g$ depends only on the image of $g$ in
the quotient group $\Gx=\mc/\Einfty$.
\end{corollary}

It is easy to see that the wreath recursion $\Phi$ induces a
well-defined wreath recursion (which we will also denote $\Phi$) on
$\Gx=\mc/\Einfty$, given by the same formula on the
generators. We also denote by $\phi$ the virtual endomorphism induced
on $\Gx$ by the virtual endomorphism $\phi$ of $G$.

We have shown earlier that the iterated monodromy group of
$F(w)=\left(\frac{2-w}{w}\right)^2$ is the quotient of $\mc$ given
by the presentation
\begin{equation}\label{eq:Gpres}\img{F}=G=\Z[i]\rtimes\Z/4=\langle
A,B|\,A^2,(AB)^4,B^4\rangle.\end{equation} The virtual
endomorphism $\phi$ induces a virtual endomorphism of $G$ (still
denoted $\phi$), which is contracting on $G$.

In $\mc$ we have the equalities
\begin{equation}
  \label{eq:firstrelGx}
  \phi(a^2)=1,\quad\phi((ab)^4)=b^{-1}a^{-2}b,\quad
\phi(a^{-1}(ab)^4a)=a^{-2},
\end{equation}
which imply that $a^2\in\mathcal{E}_1$ and $(ab)^4\in\mathcal{E}_2$.

The following lemma computes the nucleus of the wreath recursion on
$\Gx$.
\begin{lemma}
\label{lem:nucinf} For every $g\in\Gx$ there exists $n$
such that for every word $v\in X^n$ of length greater than $n$ the
restriction $g|_v$ (computed with respect to the wreath recursion
$\Phi$) belongs to the set
\[\nuke=\left\{a, ab, ab^2, ab^3, a^b, b^a, (ab)^2, b^{-2}ab,
  abab^2\right\}^{\pm 1}\cup\{b^k\}_{k\in\Z}.\]
\end{lemma}
\begin{proof}
  It is sufficient to show that $\nuke$ is state-closed, and that
  there exists $m$ such that modulo the relations $a^2=(ab)^4=1$
  obtained in~\eqref{eq:firstrelGx} we have $A|_v\subseteq\nuke$
  for all words $v$ of length $m$, where
  \[A=\nuke\cdot\{1, a, b, b^{-1}.\}\]

  Figure~\ref{fig:nuclobst} shows the Moore diagram of the set
  $\nuke$.

\begin{figure}[ht]
\centering\includegraphics{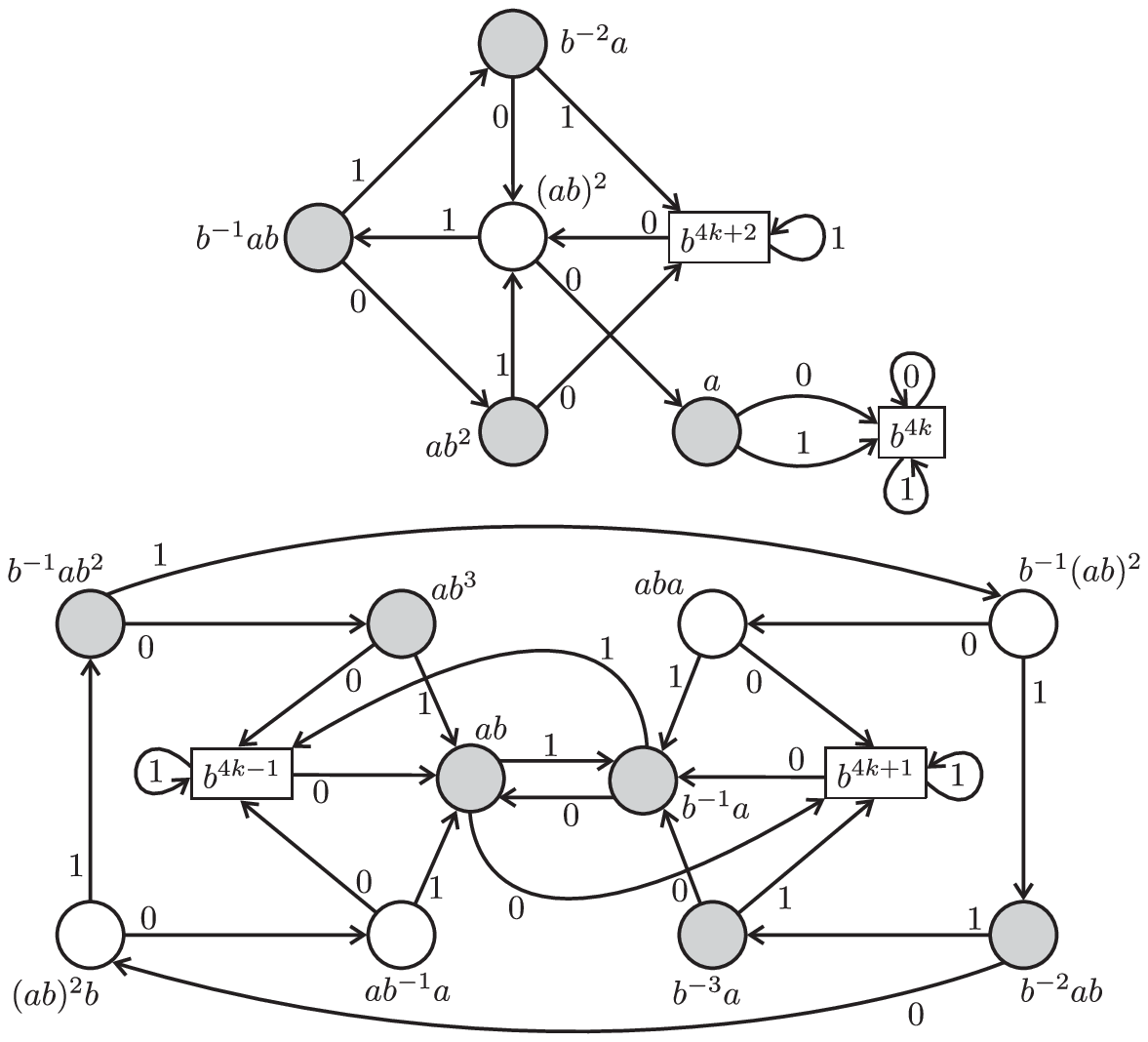} \caption{Nucleus of
$\mc$}\label{fig:nuclobst}
\end{figure}

We see that $\nuke$
  is state-closed. Therefore it is sufficient to check our
  condition for $A=\nuke\cdot\{a, b, b^{-1}\}\setminus\nuke$.
  Moreover, multiplication by $a$ from either side does not change
  the restrictions, therefore, we may take $A$ equal to
  \begin{multline*}\nuke\cdot\{b, b^{-1}\}\setminus\nuke=\{
  ab^4, abab^{-1}, b^{-2}ab^2, abab^3, b^{-1}ab^{-1},
  b^{-2}ab^{-1},\\
  b^{-3}ab,
   b^{-3}ab^{-1}, ab^{-1}ab, ab^{-1}ab^{-1},\\ b^{-1}ab^3,
  b^{-2}ab^{-1}ab, b^{-2}ab^{-1}ab^{-1}=b^{-1}aba\}\end{multline*}
  Taking into account again that multiplication by $a$ does not
  change the restrictions, and replacing some of the elements by their inverses
  (as $\nuke=\nuke^{-1}$), we reduce our checking to the set
  \[A=\{bab^{-1}, b^{-2}ab^2, bab^3,
   b^{-1}ab^3,
  b^{-2}ab^{-1}ab\}.\] But we have
  \begin{gather*}
bab^{-1}=\pair<b^{-1}ab^{-1}, bab>\sigma=\pair<(ab)^2a, a(ab)^2>\sigma\\
b^{-2}ab^2=\pair<abab^3, b^{-3}ab^{-1}a>\sigma=\pair<a\cdot bab^3,(a\cdot bab^3)^{-1}>\sigma\\
bab^3=\pair<b^{-1}ab^3, bab>\sigma=\pair<b^{-1}ab^3, a(ab)^2>\sigma\\
b^{-1}ab^3=\pair<a\cdot b^4, b^{-1}ab>\sigma\\
b^{-2}ab^{-1}ab=\pair<(ab)^2\cdot a, b^{-2}>,
  \end{gather*}
  hence we may take $m=4$, since the longest chain staying outside $\nuke$ is
  $b^{-2}ab^2\mapsto abab^3\mapsto b^{-1}ab^3\mapsto ab^4\mapsto b^4$.
\end{proof}

\begin{lemma}
\label{l:parabolic} Let $P=\bigcap_{n\ge 1}\dom\phi^n$ be
stabilizer of the path $111\ldots$ in the group $\Gx$. For every
$g\in P$ there exists $n$ such that $\phi^n(g)\in\langle
b\rangle$.
\end{lemma}

\begin{proof} Let $g\in P$ stabilize the path $111\ldots$.  There
  exists $n$ such that $\phi^n(g)\in\nuke$, where $\nuke$ is given by
  Lemma~\ref{lem:nucinf}. The sequence $(\phi^m(g))_{m\ge n}$ for must
  pass through inactive (white) states of the Moore diagram of $\nuke$
  (Figure~\ref{fig:nuclobst}) and follow arrows labeled by $1$. The
  only such infinite paths in the Moore diagram are the loops attached
  to the powers of $b$. Consequently, $\nuke\cap P=\langle
  b\rangle$.
\end{proof}

\begin{proposition}
  \label{pr:obs1} Every branched covering $f\in\mathfrak{F}$ is
  combinatorially equivalent either to $z^2+i$, or to $z^2-i$, or to
  $f_*\cdot b^k$ for some $k\in\Z$, where $b$ is as before the Dehn
  twist about the curve $\Gamma$ encircling the points $-1+i$ and
  $-i$. The branched coverings $f_*\cdot b^k$ are obstructed, with
  obstruction $\Gamma$.
\end{proposition}

\begin{proof}
  Consider the map $\ophi:\Gx\arr\Gx$ defined by
  \begin{equation}\label{eq:ophist}\ophi:g\mapsto\begin{cases}
      \phi(g) & \text{ if $w$ belongs to the domain of }\phi,\\
      a\phi(ga) & \text{ otherwise}.\end{cases}\end{equation}

  Following the line of reasoning used for the ``rabbit'' polynomial
  (see Proposition~\ref{pr:psibar}), and using
  Corollary~\ref{cor:einfty}, we check the identity $f_*\cdot
  g=\psi(g)\cdot f_*$. It then follows that the branched coverings
  $f_*\cdot g$ and $f_*\cdot\ophi(g)$ are combinatorially equivalent
  for all $g\in\Gx$.

  Lemma~\ref{lem:nucinf} and the argument
  of Proposition~\ref{pr:opsicontracting} show that it is sufficient
  to study the dynamics of the map induced by $\ophi$ on the set
  $\nuke\cup a\nuke\subset\Gx$.

  Direct computations show that the map $\ophi$ on $\nuke\cup a\nuke$
  has fixed points $a$ and $b^k$ for $k\in\Z$ and cycles
  $ab\leftrightarrow b^{-a}$ and $a^b\mapsto b^{-2a}\mapsto abab\mapsto
  a^b$.

  This shows that for every $g\in\Gx$ there exists $n$ such that
  $\ophi^n(g)\in\{a, ab, a^b\}\cup\{b^k\}_{k\in\Z}$. Consequently,
  by Corollary~\ref{cor:einfty}, every element $f\in\mathfrak{F}$ is
  combinatorially equivalent either to $f_*\cdot a=f_i\cdot a^2$, or
  to $f_*\cdot ab=f_i\cdot a^2b$, or to $f_*\cdot b^{-1}ab=f_i\cdot
  ab^{-1}ab$, or to $f_*\cdot b^k$ for some $k\in\Z$.

The branched coverings $f_i\cdot a^2$, $f_i\cdot a^2b$ and
$f_i\cdot ab^{-1}ab$ are combinatorially equivalent to,
respectively, $f_i, f_{-i}$ and $f_i$, by Theorem~\ref{th:mod5}.
The branched coverings $f_*\cdot b^k=f_i\cdot ab^k$ are
obstructed, also by Theorem~\ref{th:mod5}.
\end{proof}

\noindent The following proposition completely describes the group
$\Gx=\mc/\Einfty$.
\begin{proposition}
  $\Gx$ is a (non-split) extension of $G$ by the $G$-module
  $K=\Z[\langle B\rangle\backslash G]$. It may be given by the
  presentation
  \[\Gx=\langle a,b|\,a^2,(ab)^4,[b^4,b^{4w}]\text{ for all
  }w\in\Gx\rangle.\]
\end{proposition}

We will only need the fact that the order of $b$ in $\Gx$ is
infinite, so we give only a sketch of the proof.

\begin{proof}
  Let $N$ denote the normal closure in $\mc$ of $\langle
  a^2,(ab)^4,[b^4,b^{4w}]\text{ for all }w\in\mc\rangle$. We wish to
  prove $N=\Einfty$.

  By comparing $\mc/N$ with the presentation~\eqref{eq:Gpres}, we see
  that the map $a\mapsto A, b\mapsto B$ defines a natural epimorphism
  $\mc/N\arr G$, whose kernel $K=\langle b^{4w}:w\in\mc\rangle N$ is
  abelian, isomorphic as a $G$-module to $\Z[\langle
  B\rangle\backslash G]$.

  A straightforward argument shows that the wreath recursion
  $\Phi$ on $\mc$ induces a wreath recursion (still denoted $\Phi$) on
  the group $\mc/N$.  Furthermore, the virtual endomorphism
  $\phi:\langle B,B^A\rangle\arr G$ is a bijection, and induces a
  bijection $\Phi:K\arr K\times K$, given by
  \begin{equation}\label{eq:phibw}\Phi(b^{4w})=\begin{cases} \left(1,
        b^{4w_1}\right) & \text{if }\Phi(w)=(w_0, w_1),\\
      \left(b^{4w_1}, 1\right) & \text{if }\Phi(w)=(w_0, w_1)\sigma.
    \end{cases}
  \end{equation}

  We first prove $N\subseteq\Einfty$.  We already checked
  in~\eqref{eq:firstrelGx} that the first two generators of $N$ lie in
  $\Einfty$.  We may view Equation~\ref{eq:phibw} as a relation in
  $\Gx$; Lemma~\ref{lem:nucinf} then implies that for every $w\in\Gx$
  there exists $n_0(w)\in\N$ such that if $n\ge n_0(w)$ then all
  coordinates of $\Phi^n(b^{4w})\in \Gx^{2^n}$ are trivial, except for
  one, which is equal to $b^4$.  Therefore, the elements $b^{4w_1}$
  and $b^{4w_2}$ commute for all $w_1, w_2\in\Gx$.

  We next prove $\Einfty\subseteq N$. Assume for contradiction that
  $\mathcal E_n\supsetneq N$ for some minimal $n$; and choose
  $r\in\mathcal E_n\setminus N$. Then clearly $r\in K$, because the
  action of $r$ on $X^*$ is trivial, so the image of $r$ in $G$ is
  trivial. We have $\Phi(r)=(r_0,r_1)$, with $r_0,r_1\in\mathcal
  E_{n-1}$.  Since $\Phi$ is a bijection $K\arr K\times K$, one of
  $r_0,r_1$ is non-trivial in $\mc/N$, and we have contradicted the
  minimality of $n$.
\end{proof}

\begin{proposition}\label{pr:obs2}
  None of the maps $f_*\cdot b^k$, for $k\in\Z$, are equivalent to each other.
\end{proposition}
\begin{proof}
We may work inside the group $\Gx$, by Corollary~\ref{cor:einfty}.

  Suppose that $f_*\cdot b^r$ and $f_*\cdot b^s$ are equivalent.  Then
  there exists $h\in\Gx$ such that $h\cdot f_*\cdot
  b^r=f_*\cdot b^sh$ in $\mathfrak{F}$. We get $hb^r\cdot
  f_*=f_*\cdot b^sh$, which implies $hb^r=\phi(b^sh)$, or
  \begin{equation}\label{eq:hphi}b^{-s}hb^r=\phi(h),
  \end{equation}
  by definition of $\phi$ and
  since the left $\mc$-action on the bimodule $\mathfrak{F}$ is free
  (which follows, for instance, from the fact that the $\mc$-bimodule $\mathfrak{F}$
  is isomorphic to the bimodule associated with the self-covering
  $F(z)=\left(\frac{2-z}{z}\right)^2$ of moduli space).

By induction we get from~\eqref{eq:hphi} that for every $n$
\[\phi^n(h)=b^{-sn}hb^{rn}.\]
By Lemma~\ref{l:parabolic}, there exists $n$ such that
$\phi^n(h)\in\langle b\rangle$. Consequently, $h=b^m$ for some
$m\in\Z$, hence $b^{-s}hb^r=h$ and $b^{r-s}=1$. But this implies
that $r=s$, since the order of $b$ is infinite.
\end{proof}

\begin{corollary}
\label{cor:mod5}
  Let $f$ be a degree-two topological polynomial with preperiod of
  length $1$ and period of length $2$. Then $f$ is combinatorially
  equivalent to precisely one of $f_i$, $f_{-i}$, or $f_*\cdot b^n$ for some
  $n\in\Z$.

  Theorem~\ref{th:mod5} describes the equivalence classes of the
  polynomials $f_i$ and $f_{-i}$. If $f_*\cdot g=f_i\cdot ag$ is an
  obstructed topological polynomial, then it is equivalent to
  $f_*\cdot b^n$ if and only if there exists $m\in\N$ such that
  $\ophi^m(g)=b^n$ in $\Gx$, where $\ophi:\Gx\arr\Gx$ is given
  by~\eqref{eq:phi} and~\eqref{eq:ophist}.
\end{corollary}

\section{Preperiod $2$, period $1$}\label{s:quater}
\subsection{Iterated monodromy groups}
We next consider the degree-two topological polynomials with
ramification graph $c_1\mapsto c_2\mapsto c_3\mapsto c_3$, where
$c_1$ is the critical value. The calculations are close to those
for the ``rabbit'' and ``airplane'' polynomials, so will be given
in a more condensed manner.

If $f(z)=z^2+c$ is a polynomial whose ramification graph has preperiod
$2$ and period $1$, then the parameter $c$ must be a root of the
polynomial $x^3+2x^2+2x+2$, i.e.\ one of approximately
\[-0.2282 + 1.1151i, \quad
  -0.2282 - 1.1151i, \quad
  -1.5437.
\]
The corresponding points of the Mandelbrot set are the landing points of
the external rays at angles $1/4, 3/4$ and $5/12$, respectively.
The last angle is particular because under doubling we have
$5/12\mapsto 5/6\mapsto 2/3\mapsto1/3\mapsto2/3$; but the dynamical rays with
angles $1/3$ and $2/3$ land at the same point of the Julia set. Let us denote the
corresponding polynomials by $f_{1/4}, f_{3/4}$ and $f_{5/12}$.

\begin{figure}[ht]
\centering\includegraphics{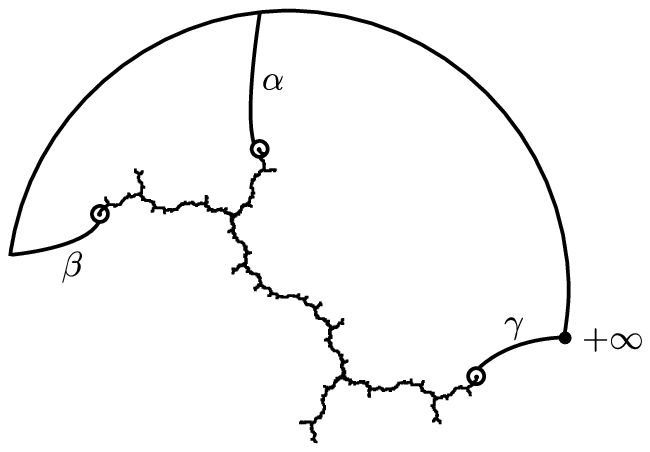}
\caption{}\label{fig:quater}
\end{figure}

The wreath recursions are given by
\begin{gather*}
\Phi_{f_{1/4}}\left(\alpha\right)=\pair<\alpha^{-1}\beta^{-1},\beta\alpha>\sigma,\quad
\Phi_{f_{1/4}}\left(\beta\right)=\pair<\alpha,1>,\quad
\Phi_{f_{1/4}}\left(\gamma\right)=\pair<\gamma,\beta>;\\
\Phi_{f_{3/4}}\left(\alpha\right)=\pair<\beta^{-1}\alpha^{-1},\alpha\beta>\sigma,\quad
\Phi_{f_{3/4}}\left(\beta\right)=\pair<1,\alpha>,\quad
\Phi_{f_{3/4}}\left(\gamma\right)=\pair<\gamma,\beta>;\\
\Phi_{f_{5/12}}\left(\alpha\right)=\pair<\alpha^{-1}\gamma^{-1},\gamma\alpha>\sigma,\quad
\Phi_{f_{5/12}}\left(\beta\right)=\pair<\alpha,1>,\quad
\Phi_{f_{5/12}}\left(\gamma\right)=\pair<\gamma^\alpha,\beta>.
\end{gather*}
Here the $\alpha$ of $f_t$ is a small loop going positively around
the critical value and connected to the circle at infinity by the
external ray at angle $t$, the $\beta$ is connected by the ray at
angle $2t$ and the $\gamma$ at angle $0$ (for $f_{1/4}$ and
$f_{3/4}$) or $2/3$ (for $f_{5/12}$). See for example the loops
$\alpha, \beta, \gamma$ for the polynomial $f_{1/4}$ on
Figure~\ref{fig:quater}.

Their respective nuclei are:
\begin{gather*}
  \nuke_{1/4} =
  \{1,\alpha,\beta,\gamma,\gamma^{\alpha\beta},\alpha\beta,\beta\alpha,(\beta\alpha\gamma)^{\pm1}\},\\
  \nuke_{3/4} =
  \{1,\alpha,\beta,\gamma,\gamma^{\beta\alpha},\alpha\beta,\beta\alpha,(\alpha\beta\gamma)^{\pm1}\},\\
  \nuke_{5/12} =
  \{1,\alpha,\beta,\gamma,\gamma^\alpha,\beta^{\gamma\alpha},\alpha\gamma,\gamma\alpha,(\beta\gamma\alpha)^{\pm1}\}.
\end{gather*}

Direct verification shows that the nuclei are different as automata.
It is not hard to prove that there is no possible obstruction, so every degree-two branched covering with
the given post-critical dynamics is equivalent to precisely one of
$f_{1/4}$, $f_{3/4}$ or $f_{5/12}$.

\subsection{The iterated monodromy group on moduli space}
As before, we force the post-critical sets to be of the form
$\{0,1,w\}$, so that $w\in\CS\setminus\{0,1,\infty\}$ represents a
point in moduli space. A polynomial with critical value $0$ and
orbit $0\mapsto1\mapsto w$ is of the form $f(z)=(az-1)^2$; then we
have
\[(a-1)^2=w_1,\quad (aw_0-1)^2=w_1,\]
so $aw_0-1=-a+1$, hence $a=\frac{2}{w_0+1}$ and
\[F(w_0)=w_1=\left(\frac{w_0-1}{w_0+1}\right)^2.\]
The fixed points of $F$ are approximately
\[-0.6478 + 1.7214i,\quad  -0.6478 - 1.7214i, \quad   0.2956.\]
The fixed point $w$ corresponds to the polynomial
$\left(\frac{2z}{w+1}-1\right)^2$, which is conjugate to
$z^2-\frac{2}{1+w}$. Direct computation shows that the above fixed
points correspond to the polynomials $f_{1/4}, f_{3/4}$ and
$f_{5/12}$, respectively.

Let us consider the first point $t\approx-0.6478 + 1.7214i$, and let
$a$ and $b$ be small loops going around the points $0$ and $1$ in the
positive direction and connected to $t$ by straight segments.  The
$F$-preimages of $t$ are $t$ itself and $1/t\approx -0.1915 -
0.5089i$. Take the connecting path $\ell_0$ to be the trivial path at
$t$, and let $\ell_1$ be a path connecting $t$ to $1/t$, intersecting
the real line only once, with the intersection point strictly between
0 and 1.

Then the iterated monodromy group of $F$ is given by the
recursions
\[\Phi(a)=\pair<1,b>\sigma,\quad\Phi(b)=\pair<b^{-1}a^{-1}, a>.\]
This wreath recursion is contracting on the free group $\langle a,
b\rangle$, and has nucleus
\[\nuke=\{1, a, b, ab, a^{-1}b\}^{\pm 1}.\]

Let $\psi$ denote the virtual endomorphism of $\mc$ corresponding
to projection on the first coordinate, i.e.
\[\psi(a^2)=b,\quad\psi(b)=b^{-1}a^{-1},\quad\psi(aba^{-1})=a.\]
Define $\opsi:\mc\arr\mc$ by
\begin{equation}\label{eq:opsiq}\opsi:g\mapsto\begin{cases}\psi(g) & \text{ if $g$ belongs to the
domain of $\psi$},\\ a\psi(ga^{-1}) & \text{
otherwise.}\end{cases}\end{equation}
\begin{lemma}
  For every $g\in\mc$, there is $n\in\N$ such that
  $\opsi^n(g)\in\{1,a,a^{-1}b\}$.
\end{lemma}
\begin{proof}
  It suffices to compute the orbits of $\opsi$ on $\nuke\cup
  a\nuke$. There are fixed points $1$ and $a$, and a cycle
  $ab^{-1}a\leftrightarrow a^{-1}b$. All the other elements of $\nuke\cup
  a\nuke$ are attracted to these cycles.
\end{proof}

\begin{figure}[ht]
\centering\includegraphics{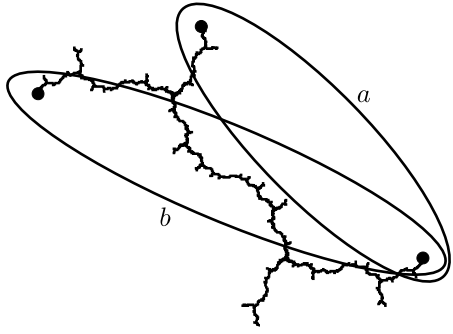} \caption{Twisting
$f_{1/4}$}\label{fig:quatert}
\end{figure}

The generators $a$ and $b$ of the mapping class group
$\mc=\pi_1(\CS\setminus\{0,1,\infty\})$ correspond to the right Dehn
twists shown in Figure~\ref{fig:quatert}, since they correspond to
the fixed post-critical point $c_3$ moving in the positive
direction around the critical value $c_1$ and the other
post-critical point $c_2$, respectively.

Computation of the iterated monodromy groups of the
maps $f_{1/4}\cdot a$ and $f_{1/4}\cdot a^{-1}b$ and their nuclei
show that the first branched covering is equivalent to $f_{5/12}$ and
the second to $f_{3/4}$.

\begin{theorem}
  Let $g\in\langle a, b\rangle$ be an arbitrary element of the mapping
  class group. Then $f_{1/4}\cdot g$ is equivalent to $f_{1/4}$,
  $f_{3/4}$ or $f_{5/12}$ if and only if the orbit of $g$ under
  iteration of $\opsi$ (given by~\eqref{eq:opsiq}) eventually reaches
  the fixed point $1$, $a$, or the cycle $ab^{-1}a\leftrightarrow
  a^{-1}b$, respectively.
\end{theorem}

\textit{Acknowledgment:} the authors are grateful to the referee for
valuable comments and suggestions.

\bibliographystyle{alpha}
\bibliography{mymath,nekrash}

\end{document}